\documentclass[fontsize=11pt,a4paper,DIV=12]{scrartcl}

\usepackage[T1]{fontenc}

\def\sc{\scshape}

\usepackage{enumerate,amsthm, amsmath, amsfonts, amssymb, mathrsfs, complexity, wasysym, graphicx,colortbl, url, hyperref, hypcap, shuffle, xargs, multicol, overpic, pdflscape, multirow, hvfloat, minibox, accents, array, multido, xifthen, xspace, ae, aecompl, blkarray, pifont, mathtools, etoolbox, dsfont, verbatim, stackrel, stmaryrd, xcolor, libertine, dynkin-diagrams, nicefrac, thmtools}

\hypersetup{colorlinks=true, citecolor=darkblue, linkcolor=darkblue,urlcolor=darkblue}
\usepackage[noabbrev,capitalise]{cleveref}

\usepackage[subrefformat=simple,labelformat=simple]{subcaption}

\definecolor{darkblue}{rgb}{0,0,0.7} 
\definecolor{darkred}{rgb}{0.7,0,0} 
\definecolor{green}{RGB}{57,181,74} 
\definecolor{violet}{RGB}{147,39,143} 
\newcommand{\darkblue}{\color{darkblue}} 
\newcommand{\red}{\color{red}} 

\newcommand{\defn}[1]{\textsl{\darkblue #1}} 

\newtheorem{theorem}{Theorem}
\newtheorem{observation}{Observation}
\newtheorem{claim}{Claim}
\newtheorem{remark}{Remark}
\newtheorem{conjecture}{Conjecture}
\newtheorem{lemma}{Lemma}
\newtheorem{proposition}{Proposition}
\newtheorem{corollary}{Corollary}

\newcommand{\fh}{facet-Hamil\-tonian\xspace}
\newcommand{\kernel}{kernel\xspace}
\newcommand{\inv}{^{-1}}
\DeclareMathOperator{\RR}{\mathbb{R}}
\DeclareMathOperator{\N}{\mathbb{N}}

\makeatletter
\def\@fnsymbol#1{%
   \ifcase#1
   \or A
   \or C
   \or F
   \or K
   \or L
   \or X
   \else
   \@ctrerr \fi
}%
\makeatother

\DeclareTOCStyleEntry[
    indent=0em,
    pagenumberformat=\normalfont,
    beforeskip=1pt plus .2pt,
    entryformat=\normalfont
    ]{tocline}{section} 
\begin{document}

\title{\huge Facet-Hamiltonicity}
\author{
Hugo Akitaya\,\thanks{University of Massachusetts Lowell, USA, hugo\_akitaya@uml.edu} \and  
Jean Cardinal\,\thanks{Universit\'e libre de Bruxelles (ULB), Belgium, jean.cardinal@ulb.be} \and  
Stefan Felsner\,\thanks{Technische Universit\"at Berlin, Germany, felsner@math.tu-berlin.de} \and 
Linda Kleist\,\thanks{Universität Potsdam, Germany, kleist@cs.uni-potsdam.de} \and
Robert Lauff\,\thanks{Technische Universit\"at Berlin, Germany, lauff@math.tu-berlin.de} \and
}

\date{}

\maketitle

\begin{abstract}
\noindent{\textbf{Abstract.}}
We consider \emph{\fh} cycles of polytopes, defined as cycles in their skeleton such that every facet is visited exactly once. 
These cycles can be understood as optimal \emph{watchman routes} that guard the facets of a polytope.
We consider the existence of such cycles for a variety of polytopes, the facets of which have a natural combinatorial interpretation.
In particular, we prove the following results:
\begin{itemize}
    \item Every permutahedron has a \fh cycle. These cycles consist of circular sequences of permutations of $n$ elements, where two successive permutations differ by a single adjacent transposition, and such that every subset of $[n]$ appears as a prefix in a contiguous subsequence. With these cycles we associate what we 
    call \emph{rhombic strips} which 
    encode interleaved Gray codes of the Boolean lattice, one Gray code for each rank. These rhombic strips correspond to simple Venn diagrams.
    \item Every generalized associahedron has a \fh cycle. This generalizes the so-called \emph{rainbow cycles} of Felsner, Kleist, Mütze, and Sering (SIDMA 2020) to associahedra of any finite type. For types $A$, $B/C$, and $D$, facets have natural interpretations in terms of arcs in triangulations, and the \fh cycles yield sequences of triangulations, where two successive triangulations differ by a single adjacent flip, and in which every arc appears and disappears exactly once. We relate the constructions to the Conway-Coxeter friezes and the bipartite belts of finite type cluster algebras.
    \item Graph associahedra of wheels, fans, and complete split graphs have \fh cycles. For associahedra of complete bipartite graphs and caterpillars, we construct \fh paths. 
    Here the facets correspond to \emph{tubes}, or connected induced subgraphs, and we obtain a sequence of elimination trees on those graphs such that every tube appears as a subtree exactly once. 
    The construction involves new insights on the combinatorics of \emph{graph tubings}.
\end{itemize}
We also consider the computational complexity of deciding whether a given polytope has a \fh cycle and show that the problem is \NP-complete, even when restricted to 
simple \hbox{3-dim}\-ensional polytopes.
\end{abstract}

\setcounter{tocdepth}{1}
\tableofcontents

\section{Introduction}
\label{sec:intro}

\begin{figure}[htb]
    \centering
    \begin{subfigure}{.3\textwidth}
        \centering\includegraphics[scale=1,page=4]{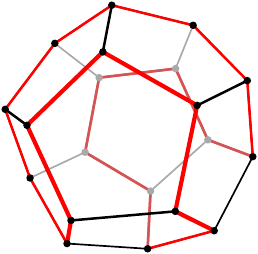}
    \subcaption{}
    \label{fig:IntroB}
    \end{subfigure}\hfil
    \begin{subfigure}{.3\textwidth}
        \centering\includegraphics[scale=1,page=5]{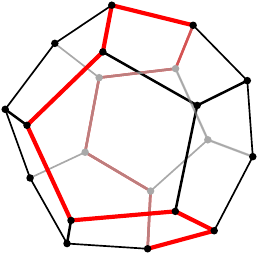}
    \subcaption{}
    \label{fig:IntroC}
    \end{subfigure}
    
  \caption{The dodecahedron with a Hamiltonian cycle ({\sc a}) and a \fh cycle ({\sc b}).}
  \label{fig:intro}
\end{figure}

Given a graph, does it contain a cycle that visits every vertex exactly once?
Such a cycle is called a \defn{Hamiltonian} cycle, in honor of Sir William Rowan Hamilton, who invented the \defn{Icosian game} in 1857. This puzzle involves finding a cycle along the edges of a dodecahedron such that every vertex is visited exactly once, see \cref{fig:IntroB}.
Hamiltonicity, the property of having a Hamiltonian cycle, has since become a fundamental theme in combinatorics and computer science, and Hamiltonicity of graphs formed by vertices and edges of polytopes, in particular, is a well-studied topic.

We propose a new notion for polytopes, that we call \defn{facet-Hamiltonicity}.
A cycle (or path) $C$ in the skeleton
of an $n$-dimensional polytope is said to be \defn{\fh} if
it visits every facet (every $(n-1)$-dimensional face) of the polytope exactly once; this means that for every facet $f$, the intersection $f\cap C$ is nonempty and connected.
\cref{fig:IntroC} illustrates a \fh cycle of the dodecahedron. A 
 polytope is \defn{\fh} if its skeleton contains a \fh cycle. 

Finding \fh cycles in simple polytopes can also be understood as a (perfect) watchman route problem where the watchman moves on the skeleton of the polytope and the polytope's surface has to be guarded. In the classical watchman route problem, one seeks a shortest closed tour of some domain such that each point of the domain is visible from some point on the tour. Clearly, a watchman route has to visit each facet at least once. 
For simple polytopes each newly visited vertex contributes one extra facet.
Measuring the length of the path by the number of edges,  a \fh cycle is a \defn{perfect} watchman route that attains this lower bound. 

\subsection{Elementary properties}

Facet-Hamiltonian cycles of simple polytopes have particularly nice properties.
Every vertex of a simple $n$-polytope $\mathcal{P}$ is incident to $n$ edges and $n$ facets.
Hence, when a path  or cycle visits a new vertex one new facet is entered, 
and one facet is left. The facet left in the next step must be different from
the one that was just entered, otherwise the cycle contains the same edge twice. It follows that the length of
a \fh cycle of a simple polytope $\mathcal{P}$ equals the 
number of facets, unless $\mathcal{P}$ has \emph{universal facets}, i.e., facets which are adjacent to all other facets. If a \fh cycle $C$ lives in the intersection of
$s$ universal facets and not more and $\mathcal{P}$ has 
$k$ facets, then the length of $C$ is $k-s$. Indeed
every facet of the simplex is universal and the $n$ simplex has \fh cycles of all length from 3 to $n+1$.
These observations are summarized below.

\begin{observation}
	\label{obs:simple}
	Let $\mathcal{P}$ be a simple $n$-dimensional polytope with $k$ facets.
	Then a \fh cycle~$C$ of $\mathcal{P}$ has the following properties:
	\begin{itemize}
        \item For every facet $f$ of $\mathcal{P}$, the intersection $f\cap C$ contains at least one edge of $\mathcal{P}$.
  \item The length of $C$ equals $k$, unless $\mathcal{P}$ 
  has $s\geq 1$ universal facets, in this case the 
  length of $C$ is between $k-s$ and $k$.
	\end{itemize}
\end{observation}

\begin{figure}[htb]
\centering
\includegraphics[height=.2\textheight]{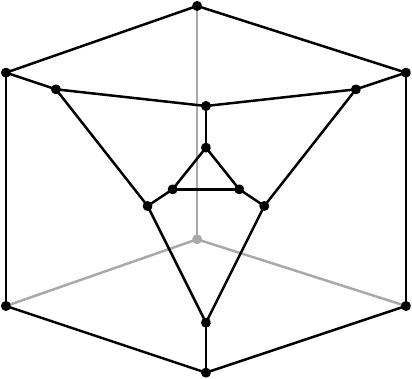}\hfil
\includegraphics[height=.2\textheight,page=2]{figures/fig-nonfhpoly2.pdf}
\caption{A three-dimensional simple polytope that is not \fh.
\label{fig:nonfh}}
\end{figure}

Unsurprisingly, not all simple polytopes are \fh. An example of non-\fh polytope is given in \cref{fig:nonfh}. To see this, let $v$ denote the gray vertex and $f$ the small triangle. Any \fh cycle uses two edges of the three cut edges highlighted in red; otherwise the facets incident to $v$ or $f$ are not visited. By symmetry, we may assume that  $e_\ell$ and $e_r$ are visited. Because they are both incident to the top facet, the \fh cycle contains either the two green or the two blue edges. This implies that either the bottom facet or $f$ is not visited.

The following question then naturally follows: What is the computational complexity of deciding whether a given polytope has a \fh cycle?
We prove that this problem is \NP-complete, even when the input polytope is three-dimensional and simple.

\begin{restatable}{theorem}{hardness}\label{thm:complexity}
The problem of deciding whether a given simple three-dimensional polytope has a \fh cycle is \NP-complete.
\end{restatable}

In what follows, we consider facet-Hamiltonicity of classical polytopes that are ubiquitous in combinatorics: \defn{permutahedra} and \defn{associahedra}.

\subsection{Permutahedra}
\label{ssec:perm}

The $(n-1)$-dimensional permutahedron~\cite{GR63} is the convex hull in $\mathbb{R}^n$ of the integer vectors denoting the permutations of $[n]$, and contained in the hyperplane of equation $\sum x_i=n(n+1)/2$. Its edges connect permutations that differ by a single adjacent transposition. 
The graph of the permutahedron is therefore the Cayley graph of the symmetric group of order $n$ for the generators consisting of adjacent transpositions.
Permutahedra are Hamiltonian, hence it is possible to list all permutations of~$[n]$ so that successive permutations differ by a single adjacent transposition.
A classical construction of such a cycle is known as the Steinhaus–Johnson–Trotter Gray code~\cite{J63,T62,T23}.

\begin{figure}[b!]
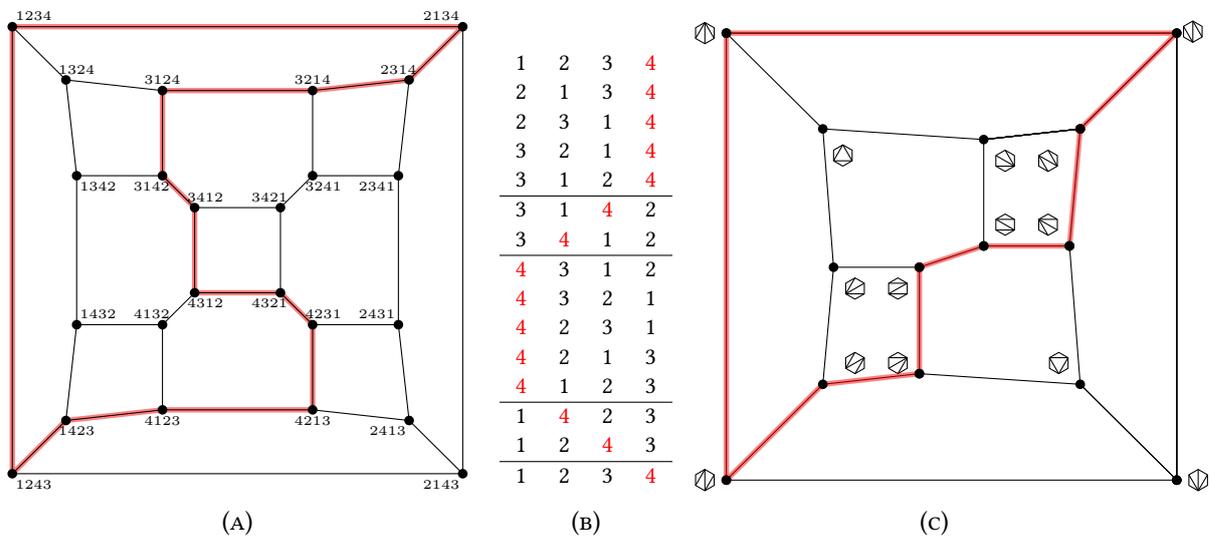

	\centering
	\begin{subfigure}{.4\textwidth}
		\centering\includegraphics[page=29]{figures/Examples.pdf}
		\subcaption{}
		\label{fig:Intro3A}
	\end{subfigure}
	\begin{subfigure}{.17\textwidth}
    \centering
	\footnotesize
    \begin{tabular}{cccc}
        1 & 2 & 3 & \red{4} \\
        2 & 1 & 3 & \red{4} \\
        2 & 3 & 1 & \red{4} \\
        3 & 2 & 1 & \red{4} \\
        3 & 1 & 2 & \red{4} \\
        \hline
        3 & 1 & \red{4} & 2 \\
        3 & \red{4} & 1 & 2 \\
        \hline
        \red{4} & 3 & 1 & 2 \\
        \red{4} & 3 & 2 & 1 \\
        \red{4} & 2 & 3 & 1 \\
        \red{4} & 2 & 1 & 3 \\
        \red{4} & 1 & 2 & 3 \\
        \hline
        1 & \red{4} & 2 & 3 \\
        1 & 2 & \red{4} & 3 \\
        \hline
        1 & 2 & 3 & \red{4} 
    \end{tabular}
		\subcaption{}
		\label{fig:Intro3C}
	\end{subfigure}
	\begin{subfigure}{.4\textwidth}
		\centering\includegraphics[page=41]{figures/Examples.pdf}
		\subcaption{}
		\label{fig:Intro3B}
	\end{subfigure}
  \caption{\label{fig:examples}
  Facet-Hamiltonian cycles in a permutahedron ({\sc a}), and an associahedron ({\sc c}). Table ({\sc b}) shows the inductive structure of the \fh cycle on the 3-dimensional permutahedron in ({\sc a}). The cycle is obtained by combining two copies of the \fh path from  $123$ to $312$ on the hexagon.}
  \label{fig:Intro3}
\end{figure}

The facets of the $(n-1)$-dimensional permutahedron are one-to-one with proper and nonempty subsets of $[n]$.
Since permutahedra are simple and no facet is incident to all other facets, Observation~\ref{obs:simple} applies, and a \fh cycle in a permutahedron must have length exactly $2^n-2$.
A \fh cycle in a permutahedron is therefore a cyclic list of $2^n-2$ permutations of $[n]$, each differing by a single adjacent transposition from its predecessor, and such that every proper subset~$S$ of $[n]$ appears as a new prefix exactly once. 
We construct \fh cycles for permutahedra. A cycle for the case $n=4$ is illustrated in \cref{fig:Intro3A}.
\begin{restatable}{theorem}{perm}\label{thm:perm}
	The $(n-1)$-dimensional permutahedron has a \fh cycle for all $n\geq 3$.
\end{restatable} 

We give a brief outline of the proof of Theorem~\ref{thm:perm}. We first prove, for every $n\geq 1$, the existence of \fh \textit{paths} between the identity permutation $1,2,\ldots ,n$ to the permutation $n,1,2,\ldots ,n-1$. This is achieved inductively, by using two copies, one of which is reversed, of the \fh paths obtained for $n-1$. The element $n$ is appended to very permutation of the first copy, and put in first position in every permutation of the second copy. We can then move from the last permutation obtained from the first copy to the first permutation in the second by moving $n$ to the front. A similar operation in the other direction closes the cycle. 
Details are given in Subsection~\ref{ssec:A-perm}. 
Figure~\ref{fig:examples} shows the construction for $n=4$.

\begin{figure}[htb]
    \centering
    \includegraphics[scale = 0.7]{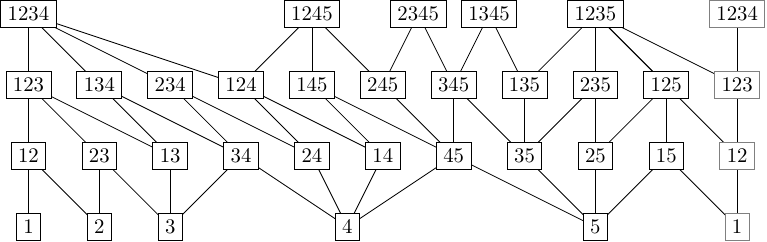}
    \caption{A cylindrically closed rhombic strip encoding \fh cycles on the 
    permutahedron for $n=5$. 
     \label{fig:rhombic3A}}
\end{figure}

\paragraph{Facet-Hamiltonian cycles and rhombic strips.}
\label{sec:rs}
Permutations are in bijection to maximal chains of the Boolean lattice. Two
permutations are adjacent on the permutahedron if the corresponding maximal
chains differ in exactly one element. Thus, a path on the permutahedron can be
encoded by a maximal chain for the first vertex and by adding a diamond detour
(a rhombic cell) for every subsequent vertex. We call such an encoding of a path a \defn{rhombic strip}. Facet-Hamiltonian cycles of the permutahedron correspond to cylindrically closed rhombic strips which are subdiagrams
of the Boolean lattice. \cref{fig:rhombic3A} shows the rhombic strip
corresponding to a \fh cycle for $n=5$. Note that the strip encodes
several cycles, for example between~$45231$ and $54321$ the cycle
may visit either $45321$ or $54231$.

More precisely, a rhombic strip in a graded poset $P$ of height $n+1$ is a 
spanning subgraph of the diagram of $P$ that admits a plane drawing with the following properties:
\begin{itemize}
    \item The vertices are placed on $n+1$ horizontal lines, labeled from $0$ to $n$ from bottom to top, so that the vertices of rank $i$ are placed on the line of label $i$, and 
    \item bounded faces are quadrilaterals.
\end{itemize}

In what follows, we will be mainly interested in rhombic strips in the Boolean lattice of subsets of~$[n]$.
Note that in our figures we omit the
vertices of rank 0 and rank $n$.
A face with two vertices of rank $k$ is
a \defn{rhombus} of rank $k$. A rhombic strip is \defn{cylindrically closed} if
on the outer face the left and the right path from rank 0 to rank $n$
are identical (here we allow duplicated vertices). Alternatively a closed
rhombic strip can be seen as a drawing on the sphere where the vertices of
rank~0 and rank~$n$ are at the south and north pole respectively, the
other ranks are represented by circles of latitude, and all the faces are
rhombi. In this paper we present rhombic strips in Figures~\ref{fig:rhombic3A},
\ref{fig:rhombic3}, \ref{fig:rhombic4}, \ref{fig:frieze}, \ref{fig:triangD}, 
\ref{fig:cycle rhombic representation} and~\ref{fig:k33rhombic}. 
With exception of the one
shown in Figure~\ref{fig:frieze} they are all closed. Furthermore, in most of these cases the ranks represent the dimension in the face lattice of an associated polytope. In the case of the type $A$ and $B$-permutahedra, these polytopes are the simplex and the cube respectively. In case of the graph associahedra (\Cref{sec:grassoc}), the rank $k$ consists of the tubes containing $k$ vertices. 

Many combinatorial polytopes can be realized as \defn{generalized permutahedra}~\cite{PRW08,P09}.
A $(n-1)$-dimensional generalized permutahedron is obtained from the $(n-1)$-dimensional permutahedron by translations of facet-defining hyperplanes. Hence the facets of
the generalized permutahedra correspond to proper subsets of $[n]$. Therefore the poset associated to the rhombic strips is a sub-poset of the Boolean lattice. Vertices of
generalized permutahedra which are incident to facets associated with subsets of all
cardinalities are called \defn{regular}. They can be encoded by a permutation or a
maximal chain in the Boolean lattice. We can therefore associate a closed rhombic strip with a \fh cycle of a generalized permutahedron that
only uses regular vertices. 

Given the rhombic strip corresponding to a \fh cycle of a generalized permutahedron $P$
we can look at the elements of rank $k$. There we find a list of all $k$-sets
which define facets. Consecutive $k$-sets have a symmetric difference of size 2.
In other words we have a Gray code for the $k$-facets of $P$. Since $k$ is arbitrary
we have a Gray code on each rank, and these Gray codes are interleaved
by the planarity condition.

In what follows, these regular cycles will play an important role. Many of our constructed cycles are of this stronger form. Especially in \Cref{sec:grassoc} some constructions only use 
nested tubings (see Section~\ref{sec:tools} for a definition), which are exactly the regular vertices.

\paragraph{Rhombic strips and Venn diagrams.}
\label{subsubsec:venn}

A \defn{Venn diagram} of the subsets of $[n]$ is a collection of $n$ simple closed curves $\gamma_1,\ldots,\gamma_n$ in the plane such that for every subset $S\subset[n]$ the region 
\[
\bigcap_{i\in S}\text{int}(\gamma_i)\cap
\bigcap_{i\in[n]\backslash S}\text{ext}(\gamma_i)
\] 
is nonempty and connected. This implies that the curves cuts the plane into regions and every region corresponds to a subset of $[n]$.
A Venn diagram is \defn{simple} if no three curves intersect in a point. See~\cite{Venn-diagrams-H63, Venn-diagrams-G75} for history and the survey~\cite{Venn-diagrams-R01} by Frank Ruskey and Mark Weston for details. 

\begin{figure}[htb]
    \centering
    \includegraphics[width=0.7\linewidth]{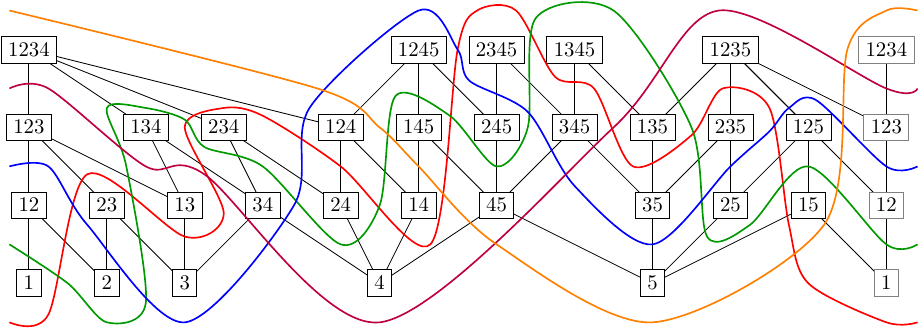}
    \caption{A rhombic strip of the 5 dimensional permutahedron and its corresponding Venn diagram. Note that the left and right side are identified as we are on the sphere.
    }
    \label{fig:venn}
\end{figure}

Considering a closed rhombic strip of the permutahedron of dimension $n$, we see all subsets of $[n]$  as the vertices of a plane graph drawn on the sphere. This graph is the dual of a corresponding Venn diagram obtained as follows. For $i\in[n]$, draw a curve with exactly those sets containing $i$ above it. To see that this is possible we draw the curve one segment at a time from left to right. Say our current endpoint is in some rhombic face. Then the bottom of this rhombus represents a set not containing $i$ while the top does contain it. Depending on whether the right vertex does contain $i$ or not we know which edge to cross next. See \Cref{fig:venn} for an illustration.

In the other direction, any simple Venn diagram has a dual graph which is a rhombic strip. Hence the rhombic strips of the permutahedron are in bijection to the simple Venn diagrams on the subsets of $[n]$. See \cite{Venn-diagrams-KCRSW04} and the figures therein for more details about simple symmetric Venn diagrams and their connections to rhombic strips. It is an open problem whether simple symmetric (with a $n$-fold rotational symmetry) Venn diagrams exist for all prime numbers $n$. Such Venn diagrams would correspond to rhombic-strips which are build from one pattern repeated $n$ times from left to right. 

\subsection{Associahedra}
\label{ssec:assoc}

Associahedra appear in various areas of mathematics~\cite{T51,S63,S63b,K01,MR1887642,MR2004457,CSZ15,AA17,PSZ23}.
They are known for encoding triangulations of a convex polygon~\cite{STT88,L89,HN99,LRS10,LP18}.
In fact, their graphs are the \defn{flip graphs} on triangulations of a convex polygon, where a flip consists of replacing the diagonal of a quadrilateral by the other diagonal.
Associahedra are known to be Hamiltonian, hence we can list all triangulations of a convex polygon so that two successive triangulations differ only by a single flip~\cite{LRR93,HN99}.
The facets of the $(n-1)$-dimensional associahedron are one-to-one with the diagonals of a convex $(n+2)$-gon.
The number of facets of the $(n-1)$-dimensional associahedron is $(n+2)(n-1)/2$.
From Observation~\ref{obs:simple}, this is the length of any \fh cycle.
Finding such a cycle amounts to finding a cyclic list of $(n+2)(n-1)/2$ triangulations of a convex $(n+2)$-gon, each differing by a single flip from its predecessor, and such that every diagonal of the polygon appears as a new edge of the triangulation exactly once.

In a series of seminal papers~\cite{MR1887642,MR2004457,MR2110627,MR2295199}, Fomin and Zelevinsky developed the theory of \defn{cluster algebras}, a field with numerous connections to other areas of mathematics. For cluster algebras of finite type, the so-called \defn{cluster complex} is the dual of a \defn{generalized associahedron} of that type, and facets of the associahedron have a natural interpretation as cluster variables. Finding \fh cycles or paths on associahedra can therefore be interpreted as generating cluster variables of a finite type cluster algebra efficiently.
The first geometric realization of generalized associahedra is due to Chapoton et al.~\cite{CFZ02}. Further generalizations to Coxeter groups are known~\cite{R06,HLT11}. 

The classical associahedra are generalized associahedra for Weyl groups of type $A$, and are therefore referred to as \defn{type $A$} associahedra.
Associahedra of type $B/C$ are also  known as \defn{cyclohedra}, 
or \defn{Bott-Taubes polytope}~\cite{S03,BT94}. 
Cyclohedra have a convenient combinatorial model defined in terms of centrally symmetric triangulations of a convex $2n$-gon, 
in which edges are symmetric pairs of diagonal flips.
Type $D$ associahedra, finally, also have nice combinatorial models, involving symmetric pseudotriangulations 
of a convex polygon minus a disk, as introduced by Ceballos and Pilaud~\cite{CP16}, or triangulations of a punctured polygon.
We show that all generalized associahedra admit a \fh cycle.
Note that \fh cycles for associahedra of type $A$ already appeared in the work of Felsner et al.~\cite{FKMS18,FKMS20} as so-called \defn{rainbow cycles}.

\begin{restatable}{theorem}{assoc}\label{thm:fhgenassoc}
  Generalized associahedra of all finite types are \fh.
\end{restatable}

\begin{figure}[htb]
	\centering
        \includegraphics[page=48]{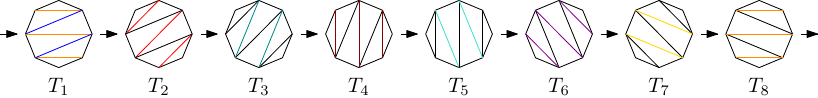}
	\caption{A \fh cycle on the 5-dimensional associahedron. Intermediate steps between the depicted triangulations are omitted, and consist of flipping triples or pairs of parallel diagonals, in any order.}
	\label{fig:associahedron}
\end{figure}

Figure~\ref{fig:associahedron} illustrates a simple construction of a \fh cycle for the type $A$ associahedron, in which two subsets of parallel diagonals are flipped alternatingly.

Permutahedra and associahedra are closely related families of polytopes, with many common generalizations.
As mentioned above, the theory of generalized permutahedra (also referred to as \defn{deformed permutahedra}) pioneered by Postnikov~\cite{PRW08,P09}, deals exactly with the properties and applications of deformations of permutahedra, which include associahedra.
We will focus on a class of generalized permutahedra called \defn{graph associahedra}, studied in particular by Carr and Devadoss~\cite{CD06,D09}.

\subsection{Graph associahedra}
Graph associahedra are polytopes defined from a given graph  that generalize permutahedra and associahedra. We first introduce the relevant terminology.

\begin{figure}[htb]
	\centering
	\includegraphics[page=43]{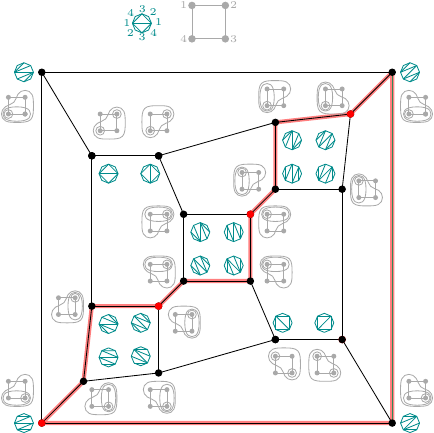}
	\caption{\label{fig:cyclo3}Illustration of $\mathcal A(C_4)$ and a \fh cycle.}
\end{figure}

Let $G=(V,E)$ be a graph. 
A \defn{tube} $t$ of $G$ is a nonempty strict subset of $V$ such that the induced subgraph $G[t]$ is connected. Two tubes $t_1\ne t_2$ are said to be
(i) \defn{non-adjacent} if $G[t_1\cup t_2]$ is not connected, and (ii) \defn{nested} if either $t_1\subset t_2$ or $t_2\subset t_1$. We call two tubes \defn{compatible} if they are either nested or non-adjacent. 
\cref{fig:compatible} illustrates examples of non-compatible and compatible tubes.
A~\defn{tubing} of $G$ is a collection of pairwise compatible tubes.
The graph associahedron $\mathcal{A}(G)$ of $G$ is a polytope whose face lattice is the inclusion poset of the tubings of $G$.
Its vertices are one-to-one with inclusion maximal tubings of $G$, and its facets are one-to-one with tubes of $G$. An example is given in \cref{fig:cyclo3}.
Since $\mathcal{A}(G)$ is simple and its facets correspond to tubes of $G$, \cref{obs:simple} implies that the length of a \fh cycle is the number of tubes of $G$, see also \cref{tab:stats}. 

\begin{table}[htb]
\caption{\label{tab:stats}Classes of graph associahedra and their number of vertices and facets}
\renewcommand{\arraystretch}{1.25}
\begin{center}
\begin{tabular}{cccc}
graph associahedra & underlying graphs & \# vertices & \# facets\\
\hline
permutahedra & $K_n$ & $n!$& $2^n-2$ 
\\[1pt]
associahedra & $P_n$ &$\frac{1}{n+1}\binom{2n}{n}$ & $(n+2)(n-1)/2$\\[3pt]
cyclohedra & $C_n$ & $\binom{2n- 2}{n-1}$& $n(n-1)$ \\[3pt]
stellohedra & $S_n$ & $\sum_{i=0}^{n-1}{i!\binom{n-1}{i}}$& $2^{n-1}+n-2$ 
\end{tabular}
\end{center}
\end{table}

The vertices of $\mathcal{A}(G)$, hence maximal tubings of $G$, can also be seen to be one-to-one with \defn{elimination trees} of $G$.
An elimination tree of $G$ is a rooted tree on the same set of vertices $V$ as $G$,
obtained by choosing any vertex $r\in V$ as root, and attaching as subtrees
the elimination trees obtained by recursing on the connected components of $G[V\setminus\{r\}]$.
The tubes of the maximal tubing are the subsets of vertices 
belonging to the proper subtrees of the elimination tree.
Edges of graph associahedra correspond to \defn{rotations} between elimination trees, or, equivalently, to \defn{flips} between pairs of tubes~\cite{CLP18,BCIKL23,CMM22}.
A \fh cycle in the graph associahedron $\mathcal{A}(G)$ is then simply a cyclic list of elimination trees of $G$, each differing by a single rotation from its predecessor, 
and such that every tube of $G$ appears as a new subtree of the elimination tree exactly once.

The $(n-1)$-dimensional permutahedron is the associahedron $\mathcal{A}(K_n)$ of the complete graph $K_n$ on~$n$ vertices.
Elimination trees in the complete graph are one-to-one with permutations of the vertices, and rotations are just adjacent transpositions.
Similarly, the classical $(n-1)$-dimensional associahedron is the associahedron $\mathcal{A}(P_n)$ of the path $P_n$ on $n$ vertices.
In that case, the elimination trees are the duals of the triangulations of a convex $(n+2)$-gon, via a classical Catalan bijection, and the rotations are the classical binary tree rotations.
The associahedron $\mathcal{A}(C_n)$ of the $n$-vertex cycle $C_n$ is the cyclohedron, and also the associahedron of type $B$.
For the $n$-vertex star $S_n$, the associahedron~$\mathcal{A}(S_n)$ is known as the \defn{stellohedron}, whose vertices are one-to-one with partial permutations (ordered subsets) of $[n]$. In that case, the rotations or flips correspond to either adjacent transpositions in the partial permutation, or adding or removing an element at the end of the partial permutation~\cite{CMM22}. 
Besides being fundamental objects in algebraic and geometric combinatorics~\cite{AA17}, graph associahedra have found applications in data structures~\cite{BCIKL23,BK22,B22} and causal inference~\cite{MUWY18,SWU21}.

Graph associahedra were shown to be Hamiltonian by Manneville and Pilaud~\cite{MP15}.
Hence it is possible to list all elimination trees of a graph in such a way that every tree differs by a single rotation from its predecessor in the list. Efficient algorithms for the case of chordal graphs were given by Cardinal, Merino, and M\"utze~\cite{CMM22}.
We show that facet-Hamiltonicity also holds for several families of graph associahedra.

\begin{restatable}{theorem}{graphAssociahedra}
  \label{thm:grassoc}
  Graph associahedra of complete graphs, paths, cycles, stars, wheels, fans, and complete split graphs  are \fh.
\end{restatable}

We give different methods for constructing \fh cycles on those polytopes, and prove
a number of interesting properties of those cycles along the way.
For associahedra of complete bipartite graphs and caterpillars, we give constructions of \defn{\fh paths}, which are defined analogously.

\begin{restatable}{theorem}{BipCatPath}\label{thm:bipcat}
    Graph associahedra of complete bipartite graphs and caterpillars have a \fh path.
\end{restatable}

Given the last result, one might wonder about the relation of \fh cycles and paths. 
In particular, is it always true that if a polytope is \fh, then it also admits a \fh path?
We show that this is indeed true for all simple 3-polytopes. However, perhaps surprisingly, this is not the case for non-simple 3-polytopes.

\begin{restatable}{theorem}{pathCycle}\label{thm:pathCycle}
    \begin{itemize}
    \item[(i)] If a simple $3$-dimensional polytope $\mathcal P$ has a \fh cycle, then it has a \fh path.
    \item[(ii)] There exists a (non-simple) 3-polytope $\mathcal P$  which has a \fh cycle but no \fh path.
    \end{itemize}
\end{restatable}

\subsection{Related work}
In this subsection we comment on connections between facet-Hamiltonicity and
various other topics in combinatorics and geometry.

\paragraph{Hamiltonian cycles.}

Studies on Hamiltonicity properties of graphs of polytopes are intimately tied to the theory of planar graphs, since from Steinitz Theorem every 3-connected simple planar graph is the graph of a 3-polytope.
Classical results in the field include Tutte's Theorem on the Hamiltonicity of 4-connected planar graphs~\cite{T56}. 
Barnette conjectured that simple 3-polytope whose facets have an even number of vertices are Hamiltonian~\cite{B69}. 
In higher dimension, a classical result of Naddef and Pulleybank states that every 0/1 polytope is Hamiltonian~\cite{NP84}; see also Merino and M\"utze~\cite{MM23b}.
Barnette conjectured that every simple $n$-polytope, for $n\geq 4$, is 
Hamiltonian (see \cite{H17}, Chapter~19).
Many \defn{Gray codes} for families of combinatorial objects are actually Hamiltonian cycles or Hamiltonian paths on graphs of polytopes, see for instance the recent survey on Gray codes from M\"utze~\cite{T23} and the series of papers on various Gray codes on polytopes generated via a simple greedy algorithm~\cite{CHMMM23, CMM22, HHMW22, HM21, MM23}.

\paragraph{Rainbow cycles.}

The special case of \fh cycles in associahedra has been studied previously in the guise of \defn{rainbow cycles} by Felsner, Kleist, M\"utze, and Sering~\cite{FKMS18,FKMS20}. Our results on generalized associahedra  extend this to associahedra of any finite type.
The original definition of a rainbow cycle in a flip graph is a cycle in which each type of flip operation occurs exactly once. 
In many settings there are several natural definitions of a flip type.
The considered rainbow cycle on permutations of $[n]$, for instance, is a cycle such that every pair of elements is swapped exactly once. In that case, the cycle has length~$\binom{n}{2}$ 
(instead of $2^n-2$ in our case) and does not live on the permutahedron.
Rainbow cycles for plane spanning trees and noncrossing matchings were also investigated in~\cite{FKMS18,FKMS20}.

\paragraph{Watchman routes.}

A watchman route of some domain is a closed tour of some domain such that each point on of the domain is visible from some point on the tour. Usually one is interested in shortest tours. A \fh cycle of a polytope can be understood as a tour where the watchman walks on the skeleton and sees each facet in a consecutive interval. For simple polytopes each newly visited vertex contributes one extra facet. Thus,  measuring the length of the path by the number of edges, the length of each watchman tour on the skeleton is lower bounded by the number of facets. Hence, a \fh cycle is an optimal watchman route in that sense.

While the watchman route problem is polynomially solvable in simple polygons \cite{watchmanSimple}, it is \NP-hard in polygons with holes~\cite{watchmanChin}. Mitchell~\cite{watchmanApproximation} presents 
an approximation algorithm with factor~$O(\log^2n)$ 
and shows that, unless $\NP=\P$, there does not exist an approximation with factor in~$o(\log n)$. Watchman routes have been also been studied for the exterior of polygonal regions~\cite{dumitrescu2023observation,ntafos1994external} as well as  for lines and line segments \cite{watchmanLines}.
In some variants the guards are restricted to walk on the boundary of the polygon, see for instance~\cite{icking1992two}.
Our hardness result for three-dimensional polyhedra (\cref{thm:complexity}) is completing this picture.

\paragraph{Hirsch conjecture.}

Hirsch conjectured that the diameter of a $n$-polytope with $t$ facets is at most $t-n$.
This conjecture was refuted in 2010 by Francisco Santos~\cite{S12}.
A \defn{nonrevisiting} path $P$ in the graph of a polytope has the property that the 
intersection of $P$ with any facet is either empty or a path on the skeleton of the 
facet. A polytope satisfies the \defn{nonrevisiting path property} if there 
exists a nonrevisiting path between any two vertices.
It is well-known that the nonrevisiting path property is equivalent to Hirsch's bound on 
the diameter, and it was once conjectured by Klee and Wolfe that every polytope 
satisfied the nonrevisiting path property~\cite{K65,K66,H03,S13}.
Several positive results have been proved by Barnette~\cite{B86,B90}.
Nonrevisiting cycles for graphs on surfaces have been considered by Pulapaka~\cite{P99}.
Our problem is a Hamiltonian counterpart of the nonrevisiting path property: Does there exist a nonrevisiting cycle visiting every facet?

\paragraph{Geodesics on graph associahedra.}

Several other properties of shortest paths, or \defn{geodesics}, on graph associahedra have been studied.
A natural question is to bound the \defn{diameter} of graph associahedra, defined as the length of the longest geodesic, hence the maximal distance of two vertices in the skeleton graph.
The diameter of associahedra has been precisely nailed down only recently~\cite{STT88,P14}.
Bounds on the diameter of several other classes of graph associahedra are known, including cyclohedra~\cite{P17}, tree associahedra~\cite{CLP18}, and complete split and bipartite graph associahedra~\cite{CPV22}. The complexity of the problem of finding a shortest path between two vertices of a graph associahedron has been studied recently~\cite{IKKKMNO23,CPV23}. Nonrevisiting properties of such shortest paths have been studied by Manneville and Pilaud~\cite{MP15}, and Ceballos and Pilaud~\cite{CP16}. It is known, in particular, that stellohedra do not satisfy the \defn{non-leaving face property}: There exist pairs of vertices that belong to a common facet, but between which all shortest paths leave and reenter this facet.

\paragraph{Visiting faces of other dimensions.}

It is natural to wonder whether it is possible to similarly construct $k$-face-Hamiltonian cycles in the skeleton of a simple polytope for any fixed  $k\in[d-1]$. Clearly, the case $k=0$ corresponds to a Hamilton cycle and the case $k=d-1$ to a \fh cycle. For instance, the \emph{equatorial} cycle of the 3-cube visits all edges. However, there is no hope to find such cycles for all $k$ in well-behaving simple polytopes such as associahedra and permutahedra. 
Indeed, if a simple $n$-polytope has a $k$-face-Hamiltonian cycle, then the total number of $k$-faces must be divisible by $\binom{n-1}{k}$, the number of new $k$-faces seen when a new vertex is visited. 
The 4-cube has 32 edges, not a multiple of 3, and the 3-dimensional associahedron has 21 edges, which is odd, in both cases the  condition is violated.
The divisibility condition for $k=1$ is always true for the permutahedron, but one can show that for instance the 3-dimensional permutahedron has no 1-face-Hamiltonian cycle.

\subsection{Open problems}
\label{sec:discussion}

We now discuss a number of open problems.

\paragraph{Type \texorpdfstring{$B$}{B} permutahedra.}

So far we only found rhombic strips for the type $B$ permutahedra in dimensions~3 and 4. 
We believe, however, that they exist in all dimensions.
This would imply the following.

\begin{conjecture}
Type $B$ permutahedra of all dimensions are \fh .
\end{conjecture}

\paragraph{Rhombic strips and truncated polytopes.}

We already commented in Section~\ref{sec:rs} on the correspondence between \fh cycles in type $A$ permutahedra and rhombic strips in the Boolean lattice. This correspondence actually holds in much more general contexts. 

\begin{figure}[htb]
\centering
    \includegraphics[page=2,scale=.7]{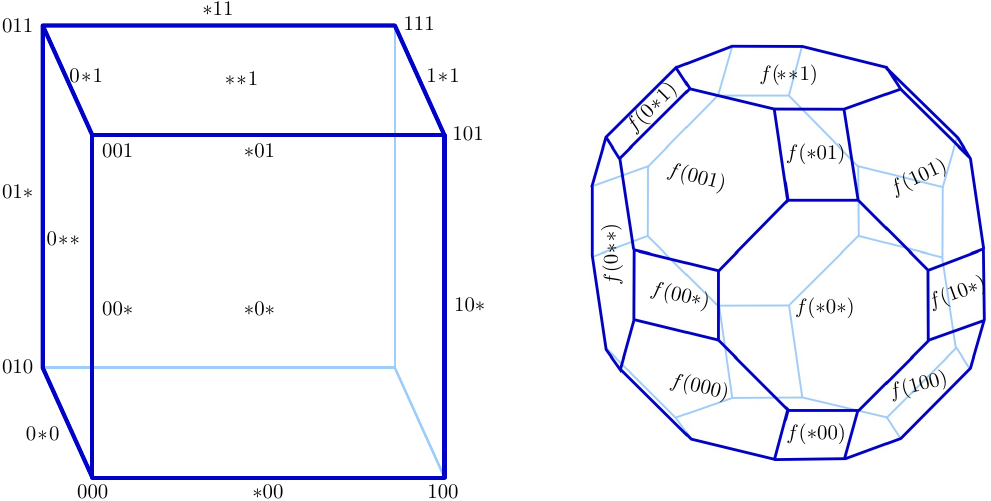}
    \caption{A 3-cube $Q$ and the 3-dimensional B-permutahedron (the truncation $Q^T$).\label{fig:trunc-cube}}
\end{figure}

Let us briefly introduce a more general construction.
Given a polytope $P$ we can ask for rhombic strips in the diagram of the face lattice of $P$. 
Such strips yield strongly restricted walks on the faces of each dimension. For example the sequences of elements in the lowest and highest rank
are Hamiltonian cycles of $P$ and its dual, respectively. 
We can also map such a rhombic strip to a
\fh cycle of the polytope $P^T$ obtained from $P$ by truncating all its proper faces. Vertices of~$P^T$ are in bijection to flags of $P$, hence maximal chains in the face lattice of $P$, and facets of~$P^T$ are in bijection to faces of $P$. 
For example, the type $B$ permutahedron can be obtained by truncating all proper faces of the hypercube, see \cref{fig:trunc-cube}. 

We believe that rhombic strips have a lot of potential for future insights and results. They tie together many results from this paper, for example by giving a connection between the bipartite belts in cluster algebra (see \Cref{sec:cluster}) and the constructions for graph associahedra (from \Cref{sec:grassoc}). In fact, most paths and cycles presented here stem from a rhombic strip.

\paragraph{Graph associahedra.}

Is it true that for every connected graph $G$, the graph associahedron $\mathcal A(G)$ is \fh?
Note that when $G$ is a cycle, we obtained a \fh cycle consisting only of \defn{nested tubings}: tubings that consist of pairwise nested tubes. Our \fh cycles for the star and the path do not have this property. This is no coincidence. 
\begin{observation}
If the graph associahedron $\mathcal A(G)$ of graph $G$ has a \fh cycle $C$ consisting only of nested tubings, then $G$ is Hamiltonian.
\end{observation}
\begin{proof}
Because all tubings of $C$ are nested and $C$ visits all facets, every vertex of $G$ is the \kernel (a tube consisting of this single vertex) of some tubing in $C$. Tracing these kernels along the \fh cycle $C$ gives a Hamiltonian cycle of $G$. This is because flipping some kernel to a different kernel happens inside a tube of size two, which guarantees an edge in~$G$.
\end{proof}

The reverse statement is not true: the Hamiltonicity of $G$ is not sufficient for the existence of a \fh cycle consisting only of nested tubings. 
Consider $G=K_4\backslash\{e\}$, the complete graph on four vertices minus an edge. 
While $G$ is Hamiltonian, $\mathcal A(G)$ does not have a \fh cycle using only nested tubings as can easily be checked by hand.

The following question follows naturally: Which graph associahedra have \fh cycles or paths consisting only of nested tubings?

\subsection{Plan of the paper}

In \cref{sec:permutahedra}, we present \fh cycles for permutahedra and thus prove \cref{thm:perm}.
We also give a construction of rhombic strips, hence of \fh cycle, for type $B$ permutahedra of dimension 3 and 4.
In \cref{sec:assoc}, we prove \cref{thm:fhgenassoc} for generalized associahedra of type $A$, $B/C$, and $D$. 
The connection between these results and known ideas in cluster algebras is developed in Section~\ref{sec:cluster}, which contains
a simple proof of the existence of \fh cycles for associahedra of all finite types. This proof makes use of a classical tool in cluster algebra, known as bipartite belts, which happen to correspond to rhombic strips for all generalized associahedra. 
In \cref{sec:grassoc}, we discuss graph associahedra. In particular, we present \fh cycles for associahedra of several graph families (\cref{thm:grassoc}) and \fh paths for complete bipartite graphs and caterpillars (\cref{thm:bipcat}). In~\cref{sec:pathsCycles}, we give the proof of \Cref{thm:pathCycle}.
Finally, we establish the \NP-completeness of deciding the existence of \fh cycles even for three-dimensional polyhedra (\cref{thm:complexity}) in \cref{sec:hardness}. 

\section{Facet-Hamiltonian cycles in permutahedra}
\label{sec:permutahedra}

We discuss  \fh cycles in permutahedra of type $A$ in \cref{ssec:A-perm} and of type $B$ in \cref{ssec:B-perm}

\subsection{Type \texorpdfstring{$A$}{} Permutahedra}
\label{ssec:A-perm}

We denote by $\mathcal A(K_n)$ the $(n-1)$-dimensional permutahedron defined by
\[
\mathcal A(K_n) = \mathrm{conv} \{(\pi(1),\pi(2),\ldots ,\pi(n)) : \pi\in S_n\},
\]
where $S_n$ is the set of permutations on $n$ elements. 
Note that the notation $\mathcal A(K_n)$ relies on the fact that the permutahedron is the associahedron of the complete graph, anticipating on Section~\ref{sec:grassoc}. 
It is well-known that the edges of the permutahedron are in bijection with pairs of permutations that differ by a single adjacent transposition. 
Moreover, the facets of $\mathcal A(K_n)$ can be labeled by subsets of~$[n]$ such that the facets incident to a vertex appear as a prefix in the permutation.

We now construct \fh cycles in permutahedra, and give a complete proof of \cref{thm:perm}. Recall the sketch and \cref{fig:Intro3A} from Subsection~\ref{ssec:perm}. As our proof is by induction, we will use the following notation to lift permutations of $[n-1]$ to permutations of $[n]$.
Let $\pi$ be a permutation of $[n-1]$, then $\pi_k$  denotes the permutation where $n$ is inserted at the $k$th position in $\pi$ for any fixed $k\in \{1,\dots, n\}$.   
Hence for $\pi=1,2,\dots, n-1$, we have $\pi_1=n,1,2,\dots, n-1$ and $\pi_n=1,2,\dots, n-1,n$.

\begin{lemma}
\label{lem:fhpaths}
 The $(n-1)$-dimensional permutahedron has a \fh path from $1, \ldots, n$ to $n, 1, \ldots, n-1$ for all $n\geq 1$.   
\end{lemma}
\begin{proof}
        We identify the facets of $\mathcal A(K_n)$
    with the nontrivial subsets of $[n]$.
	We prove the existence of the paths by induction. For $n=1$ and $n=2$, the statement is obvious. For the induction step, consider  $n\ge 3$ and let $P$ be a  \fh path in $\mathcal A(K_{n-1})$ from $\rho:=1,\dots, n-1$ to $\tau:=n-1,1,\dots, n-2$. Let $P^{n}$ denote the path obtained by replacing each permutation $\pi$ in $P$ with $\pi_n$.
	Then $P^{n}$ is a well-defined path in $\mathcal A(K_n)$ from $\rho_n$ to $\tau_n$ that introduces every facet not containing~$n$ 
 (except for the prefixes of $\rho_n$). Let $Q$ be the path $\tau_n,\tau_{n-1},\dots,\tau_1$, in which $n$ is shifted to the front. Note that $Q$ introduces all facets that are prefixes of $\tau_1$.
	 Let $\overleftarrow{P^1}$ denote the path $P$ in reverse where each permutation $\pi$ is replaced by $\pi_1$. Then,  $\overleftarrow{P^1}$ is a path in $\mathcal A(K_n)$ that introduces every facet containing~$n$ (except of the prefixes of $\tau$). The concatenation $PQ\overleftarrow{P^1}$ is a \fh path from~$\rho_n$ to $\rho_1$ as claimed.
\end{proof}

\perm*

\begin{figure}[htb]
	\centering
	\begin{subfigure}{.3\textwidth}
		\centering
		\includegraphics[page=4]{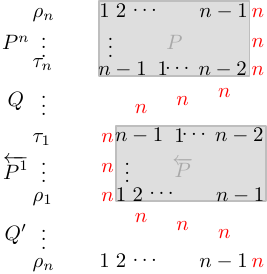}
		\subcaption{Notations}
		\label{fig:permA}
	\end{subfigure}
	\hfil
	\begin{subfigure}{.12\textwidth}
 \vbox{}\vskip-3mm\vbox{}
 \centering
     \begin{tabular}{cccc}
        \cellcolor{lightgray} 1 &\cellcolor{lightgray} 2 & \cellcolor{lightgray}  3 & \red{4} \\
        \cellcolor{lightgray} 2 &\cellcolor{lightgray}  1 &\cellcolor{lightgray}  3 & \red{4} \\
        \cellcolor{lightgray} 2 &\cellcolor{lightgray}  3 &\cellcolor{lightgray}  1 & \red{4} \\
        \cellcolor{lightgray} 3 &\cellcolor{lightgray}  2 &\cellcolor{lightgray}  1 & \red{4} \\
        \cellcolor{lightgray} 3 &\cellcolor{lightgray}  1 &\cellcolor{lightgray}  2 & \red{4} \\
        \hline
        3 & 1 & \red{4} & 2 \\
        3 & \red{4} & 1 & 2 \\
        \hline
        \red{4} &\cellcolor{lightgray}  3 &\cellcolor{lightgray}  1 &\cellcolor{lightgray}  2 \\
        \red{4} &\cellcolor{lightgray}  3 &\cellcolor{lightgray}  2 &\cellcolor{lightgray}  1 \\
        \red{4} &\cellcolor{lightgray}  2 &\cellcolor{lightgray}  3 &\cellcolor{lightgray}  1 \\
        \red{4} &\cellcolor{lightgray}  2 &\cellcolor{lightgray}  1 &\cellcolor{lightgray}  3 \\
        \red{4} &\cellcolor{lightgray}  1 &\cellcolor{lightgray}  2 &\cellcolor{lightgray}  3 \\
        \hline
        1 & \red{4} & 2 & 3 \\
        1 & 2 & \red{4} & 3 \\
        \hline
        1 & 2 & 3 & \red{4} 
    \end{tabular}
		\subcaption{$n=4$}
		\label{fig:permB}
	\end{subfigure}
	
	\caption{Construction of \fh cycles in permutahedra.}
	\label{fig:perm}
\end{figure}

\begin{proof} 
    We build the cycles from the paths defined in Lemma~\ref{lem:fhpaths}.
	To obtain a cycle in case $n\geq 3$, we add the path $Q'$ where $n$ is shifted to the back, i.e., $Q':=\rho_1,\rho_2,\ldots,\rho_n$. This introduces all facets  of $\rho_n$ and  closes the path to the \fh cycle 
	$
	\underbrace{\rho_n,\dots \tau_n}_{P^n}, \tau_{n-1}\dots,\underbrace{\tau_1,\dots, \rho_1}_{\overleftarrow{P^1}},\ldots,\rho_n,
        $
 see also \cref{fig:permA}.
	By construction, $Q'P$ visits all facets not containing $n$ and $Q\overleftarrow{P^1}$ visit all facets containing~$n$. \cref{fig:Intro3A,fig:permB} illustrate the resulting cycle for $n=4$.
\end{proof}

\subsection{Type \texorpdfstring{$B$}{B} permutahedra}
\label{ssec:B-perm}

The type $B$ permutahedron is defined as the convex hull of the points corresponding to \defn{signed permutations}, of the form $(\pm\pi(1),\ldots,\pm\pi(n))\in\RR^n$ where $\pi$ is a permutation of $[n]$. 
Equivalently, it is the \defn{zonotope} of the type $B$ root system.

\begin{figure}[hbt]
    \centering
    \includegraphics[scale = 0.7]{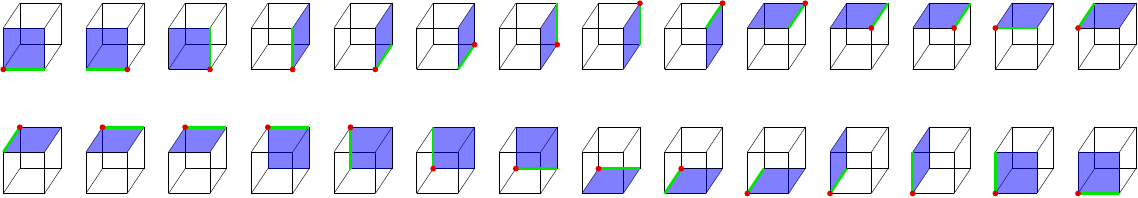}
    \caption{A \fh cycle on the three-dimensional type $B$ permutahedron represented
    by flags of the cube.}
    \label{fig:flag representation}
\end{figure}

There is a bijection between signed
permutations and maximal chains in the face lattice of the cube or
equivalently flags of the cube, see also \cref{fig:flag representation}: Faces of the $n$-cube can be encoded as
vectors of length $n$ with entries in $\{0,1,x\}$, the number of entries of
type $x$ is the dimension of the face.  Chains containing a face encoded by
$(a_1,\ldots,a_n)$ correspond to permutations with a suffix consisting of the
set $T$ defined as follows: If $a_i=1$ then $+i\in T$, if $a_i=0$ then $-i\in
T$, and if $a_i=x$ then $\pm i\not\in T$. A maximal chain prescribes all
suffices and hence a signed permutation. For example the chain
$(010),(01x),(x1x)$ corresponds to the signed permutation $(-3,-1,+2)$. 

\begin{figure}[hbt]
    \centering
    \includegraphics[scale = 0.78]{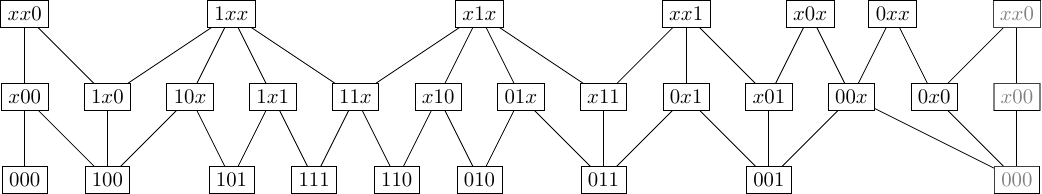}
    \caption{The \fh cycle on the type $B$ permutahedron of dimension~3 represented by its rhombic strip. There is a unique strip up to graph isomorphisms. The two cut vertices in the middle rank show that this strip encodes four \fh cycles.}
    \label{fig:rhombic3}
\end{figure}

\begin{figure}[hbt]
    \centering
    \includegraphics[scale = 0.3]{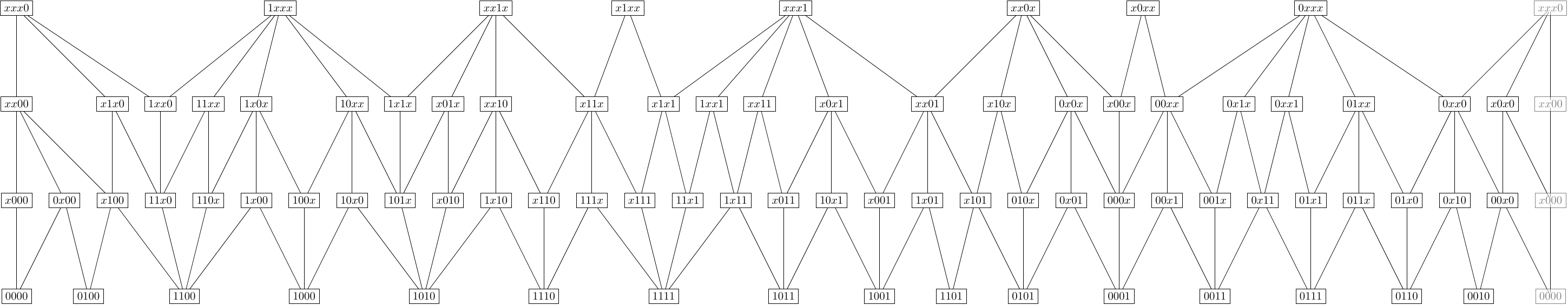}
    \caption{A rhombic strip in the type $B$ permutahedron of dimension~4.
    }
    \label{fig:rhombic4}
\end{figure}

Signed permutations are adjacent in the skeleton of the type $B$ permutahedron if
they differ either by an adjacent transposition preserving signs, or by the sign of the first element. 
The chains corresponding to adjacent signed permutations differ by exactly
one element. Hence paths on the type~$B$ permutahedron can be represented
by rhombic strips in the face lattice of the cube.
Figures~\ref{fig:rhombic3} and~\ref{fig:rhombic4} show such rhombic strips for the 3 and 4-dimensional type $B$ permutahedra. This implies the following.

\begin{proposition}
    Type $B$ permutahedra of dimension 3 and 4 are \fh .
\end{proposition}

Note that a rhombic strip for the type $B$ permutahedron yields a collection of interleaved Gray codes, one for each rank.
On rank 0 we have a standard Gray code on binary words, and on rank $k>0$
two vectors on the alphabet $0,1,x$ are adjacent if at one position an $x$ is
exchanged by $0$ or $1$ and at another position a $0$ or $1$ is made an $x$.

\section{Facet-Hamiltonian cycles in associahedra}
\label{sec:assoc}

We now present \fh cycles of associahedra of type $A$, $B/C$, and $D$, namely we prove \cref{thm:fhgenassoc} for the three main Coxeter types. The three types are considered in  Propositions~\ref{lem:asso}, \ref{lem:cyclo}, and~\ref{lem:assoD}, respectively.

\subsection{Type \texorpdfstring{$A$}{A} associahedra}
\label{sec:Aassoc}

Facet-hamiltonicity of associahedra was first proven by Felsner, Kleist, M\"utze, and Sering~\cite[Theorem 1]{FKMS18,FKMS20}. 
Below we give two proofs, the second one details the proof sketched in \cref{fig:associahedron} in Subsection~\ref{ssec:assoc}. 
Both proofs have nice interpretations in triangulations and the ideas can also be used to construct \fh cycles for associahedra of type $B/C$ and~$D$.  

We denote the $(n-1)$-dimensional associahedron by $\mathcal A(P_n)$ (again, anticipating on Section~\ref{sec:grassoc}, from the fact that the associahedron is the graph associahedron of a path on $n$ vertices).
We rely on the well-known facts that the vertices of $\mathcal A(P_n)$ are in bijection with the triangulations of a convex $(n+2)$-gon, the edges are in bijection with pairs of triangulations differing by a single flip of a diagonal, and the facets are in bijection with the diagonals.

\begin{restatable}{proposition}{asso}\label{lem:asso}
	For all $n\geq 3$, $\mathcal A(P_n)$ has a \fh cycle. 
\end{restatable}

\begin{proof}[Proof 1 -- Sketch of the proof by Felsner, Kleist, M\"utze, and Sering~\cite{FKMS18,FKMS20}]
We label the corners of the convex $(n+2)$-gon clockwise by the integers $1,2,\dots, n+2$.
Let $S_i$ denote the triangulation which contains the diagonals $\{i,k\}$ for $k\in \{i+2,\dots, n+2\}$ and $\{k,n+2\}$ for $k\in \{2,\dots, i\}$ as illustrated in \cref{fig:associahedron1A}; note that the diagonals form a bi-star. 

\begin{figure}[htb]
	\centering
	\begin{subfigure}{.15\textwidth}
		\centering
		\includegraphics[page=35]{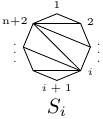}
		\subcaption{}
		\label{fig:associahedron1A}
	\end{subfigure}\hfill
	\begin{subfigure}{.37\textwidth}
		\centering
		\includegraphics[page=36]{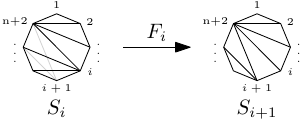}
		\subcaption{}
		\label{fig:associahedron1B}
	\end{subfigure}\hfill
	\begin{subfigure}{.37\textwidth}
		\centering
		\includegraphics[page=37]{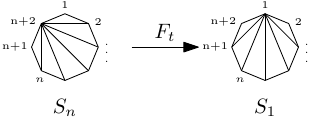}
		\subcaption{}
		\label{fig:associahedron1C}
	\end{subfigure}
	
	\caption{Illustration for proof 1 of \Cref{lem:asso} for $n=6$.}
	\label{fig:associahedron1}
\end{figure} 

There exists a simple flip sequence $F_i$ from $S_{i}$ to $S_{i+1}$ where for increasing $k\in \{2,\dots,n+1-i\}$ the diagonal $\{i,i+k\}$ is replaced by $\{i+1,i+k+1\}$, see also \cref{fig:associahedron1B}. Similarly, there exists a simple flip sequence $F_t$ from $S_{n}$, containing all diagonals of type $\{n+2,k\}$, to $S_{1}$, containing all diagonals of type $\{1,k\}$, in which for increasing $k\in \{2,\dots,n\}$ the diagonal $\{n+2,k\}$ is replaced by $\{1,k+1\}$, see also \cref{fig:associahedron1C}.
The concatenation $C$ of $F_1,F_2,\dots,F_{n-1},F_t$ is a \fh cycle. To this end, note that the diagonal $\{i,j\}$, $1 < i<j$, is not contained in $S_k$ for $k < i$ but contained in $S_i$. Hence it is introduced in~$F_{i-1}$.
Diagonals containing $1$ are introduced in 
$F_t$. Quite obviously no two triangulations 
in~$C$ coincide.
\end{proof}

We now present the idea of an alternative approach.

\begin{proof}[Proof 2]
A \fh cycle of $\mathcal A(P_n)$ has length 
$\binom{n+2}{2}-(n+2)=\frac{1}{2} (n+2)(n-1)$.
The idea is to introduce classes of parallel diagonals one after the other. For an illustration see \cref{fig:associahedron}. Note that the set of diagonals of the $(n+2)$-gon partitions into $n+2$ classes $P_i$ of parallel slopes. Two classes are \defn{compatible} if they contain no crossing diagonals and thus yield a triangulation. It is easy to see that for any class there exist exactly two other classes with which it is compatible. Thus, we may consider a cyclic list $P_1,\ldots,P_{n+2}$ such that $P_i$ and $P_{i+1}$ are compatible. We define the triangulation $T_i=P_i \cup P_{i+1}$. The triangulation $T_i$ can be transformed into $T_{i+1}$ by removing the diagonals of $P_i$ and introducing the diagonals of $P_{i+2}$. Note that the flips from $T_i$ to $T_{i+1}$ are independent and can be performed in any order.

Let $F_{i}$ denote a flip sequence from $T_i$ to $T_{i+1}$.
Then, the concatenation of $F_1,F_2,\dots, F_{n+2}$ yields a cycle $C$ because all triangulations differ. While $F_i$ and $F_{i+1}$ have different lengths if $n$ is even, $|F_i|+|F_{i+1}|=n-1$ holds in all cases because all diagonals of $T_{i+2}$ are introduced; for odd $n$, we have $|F_i|=(n-1)/2$. The cycle thus has length $(n-1)(n+2)/2$ as desired. Moreover, every inner diagonal of the $(n+2)$-gon appears exactly in one $P_i$  and is thus introduced in $F_{i-1}$. Because the cycle $C$ has a length that equals the number of diagonals, every diagonal is introduced exactly once, i.e., every facet is visited. Consequently, $C$ is a \fh cycle of $\mathcal A(P_n)$.
\end{proof}

\begin{figure}[htb]
	\centering
        \includegraphics[page=28]{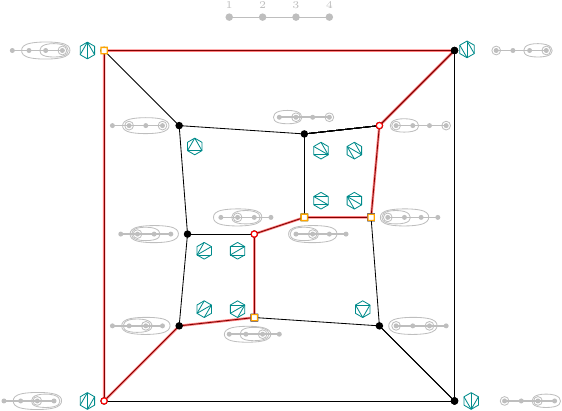}
	\caption{Illustration for a cycle that could be constructed by the proofs of \Cref{lem:asso}; square vertices are the triangulations $S_i$ in proof 1, round-hollow  vertices are the triangulations $T_i$ from proof 2.}
	\label{fig:associahedron3}
\end{figure} 

\subsection{Type \texorpdfstring{$B/C$}{B/C} associahedra}
\label{sec:Bassoc}

As mentioned before the associahedra of type $B/C$ are also known as cyclohedra and Bott-Taubes polytopes, and as the graph associahedra $\mathcal{A}(C_n)$ of the $n$-vertex cycle~$C_n$~\cite{BT94,S03}. 

We use the fact that the vertices are in bijection with the  centrally symmetric triangulations of a convex $2n$-gon, edges of the cyclohedron correspond to flipping pairs of diagonals in the triangulation or a longest diagonal, and facets correspond to pairs of symmetric diagonals. This model allows us to use ideas similar as in the proofs of \cref{lem:asso}. \cref{fig:cyclo3} illustrates the flip graph  and a \fh cycle for $\mathcal A(C_4)$.

\begin{restatable}{proposition}{cyclo}\label{lem:cyclo}
	For all $n\geq 3$, $\mathcal A(C_n)$ has a \fh cycle. 
\end{restatable}

\begin{figure}[htb]
	\centering
	\begin{subfigure}{.18\textwidth}
		\centering
		\includegraphics[page=46]{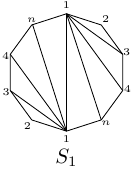}
		\subcaption{}
		\label{fig:cycle3A}
	\end{subfigure}\hfil
	\begin{subfigure}{.35\textwidth}
		\centering
		\includegraphics[page=47]{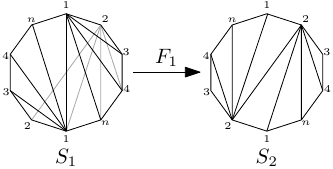}
		\subcaption{}
		\label{fig:cycle3B}
	\end{subfigure}\hfil
	\caption{Illustration for proof 1 of \cref{lem:cyclo}.}
	\label{fig:cycle3}
\end{figure}

\begin{proof}[Proof 1] Label the vertices of the $2n$-gon by $1,\dots,n,1\dots,n$ clockwise, as in \cref{fig:cycle3A}.
Let $S_i$ denote the triangulation which contains the clockwise diagonals $\{i,k\}$ for $k\in \{i+2,\dots, i+n\}$ (where $n+j=j$). Similar to the above, there is a simple flip sequence $F_i$ of length $n-1$ from~$S_i$ to~$S_{i+1}$ where the (pairs of) diagonals can be introduced by increasing length see \cref{fig:cycle3B}. Then the concatenation $F_1,F_2,\dots,F_{n}$ is \fh-cycle of length 
	$n(n-1)$: Firstly, any two triangulations are distinct as they either their longest diagonals differ, or they contain a same diagonal~$\{i,i\}$ and belong to $F_i$. Secondly, each diagonal is contained in exactly one~$S_i$.
\end{proof}

\begin{figure}[htb]
		\centering	
		\includegraphics[page=49]{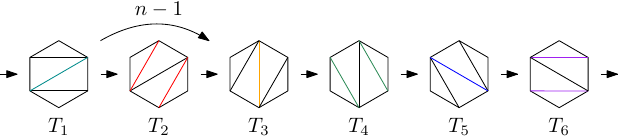}
		\caption{Illustration for proof 2 of \Cref{lem:cyclo} for $n=3$.}
		\label{fig:cyclo}
	\end{figure}

 \begin{proof}[Proof 2]
	The idea is to introduce classes of parallel diagonals one after the other in the $2n$-gon. To this end, let $P_1,\ldots,P_{2n}$ be cyclic list of classes of parallel slopes  of the $2n$-gon such that $P_i$ and $P_{i+1}$ are compatible as in the second proof of  \cref{lem:asso}.  For an illustration, consider \cref{fig:cyclo}. We define the triangulation $T_i=P_i \cup P_{i+1}$. As before, there is a flip sequence $F_i$ from  the triangulation~$T_i$ to~$T_{i+1}$ that removes the diagonals of $P_i$  and introduces the diagonals of $P_{i+2}$ in pairs; however this time, we introduce all but the longest diagonal in pairs. Thus,  $|F_i|+|F_{i+1}|=(2n-3+1)/2=n-1$.
	Then, the concatenation of $F_1,F_2,\dots, F_{2n}$ yields a \fh cycle of length $n(n-1)$. 
\end{proof}

\subsection{Type \texorpdfstring{$D$}{D} associahedra}

\begin{figure}[htb]
	\centering
	\begin{subfigure}{.2\textwidth}
		\includegraphics[page=19]{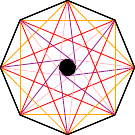}
		\subcaption{}
		\label{fig:typeDConstructionA}
	\end{subfigure}\hfil
	\begin{subfigure}{.2\textwidth}
		\includegraphics[page=40]{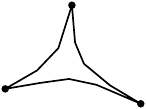}
		\subcaption{}
		\label{fig:typeDConstructionB}
	\end{subfigure}\hfil
	\begin{subfigure}{.2\textwidth}
		\includegraphics[page=38]{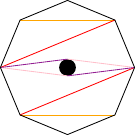}
		\subcaption{}
		\label{fig:typeDConstructionC}
	\end{subfigure}
	
	\caption{Illustration for associahedra of type $D$. ({\sc a}) The set of chords of $P$. ({\sc b}) A pseudotriangle. ({\sc c}) The zigzag pseudotriangulation $T_0$.}
	\label{fig:typeDConstruction}
\end{figure}

Type $D$ associahedra allow for nice combinatorial models where the vertices correspond to triangulations. 
Here we recall the one given by Ceballos and Pilaud~\cite{CP16}. 
(See \cref{sec:friezestofh} and Figure~\ref{fig:triangD} for a different model.)
For the $n$-dimensional associahedron of type $D$, denoted by  Asso$(D_n)$, consider the regular $2n$-gon~$P$, together with a disk $O$ placed at its center, the radius of which is small enough such that $O$ only intersects the long diagonals of~$P$. As \defn{chords} of~$P$, we consider all diagonals disjoint from $O$, together with two tangents from each vertex of $P$ to~$O$. The set of all chords is depicted in \cref{fig:typeDConstructionA}.

The vertices of  Asso$(D_n)$ are in bijection with the 
centrally symmetric \defn{pseudotriangulations} each of which contains exactly $2n$ chords. A \defn{pseudotriangulation} is a partition of a convex polygon into \defn{pseudotriangles}, defined as regions with three convex corners and an arbitrary number of reflex vertices. As in our setting reflex vertices can only come from $O$, pseudotriangles of $P$ contain at most one chain of reflex vertices, and 0,1, or 2 convex vertices on $O$.
The edges of  Asso$(D_n)$ correspond to 
flips of centrally symmetric pairs of chords and the facets correspond to centrally symmetric pairs of (internal) chords. The number of facets, and thus the length of a \fh cycle,  is $n^2$.

\begin{figure}[htb]
 	\centering
 	\includegraphics[page=20]{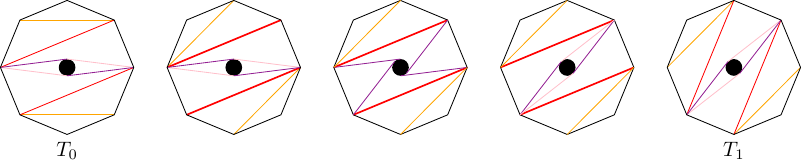}
 	\caption{Illustration for the proof of \cref{lem:assoD}: Flip sequence from $T_0$ to $T_1$ in $n$ steps for $n=4$.}
 	\label{fig:typeDProof}
 \end{figure}

While both proofs of \cref{lem:cyclo} generalize straightforwardly, we present a sketch of the second.
\cref{fig:typeDProof} depicts a flip sequence from one \defn{zigzag}-triangulation to its rotated copy for the case of $n=4$.

 \begin{restatable}{proposition}{assoD}\label{lem:assoD}
 	For all $n\geq 3$, Asso$(D_n)$ has a \fh cycle. 
 \end{restatable}

 \begin{proof}[Proof sketch]
 We consider the zigzag triangulation $T_0$ of $P$, where the long diagonal is replaced by the four tangents as illustrated in \cref{fig:typeDConstructionC}. Let $T_i$ be obtained from $T_0$ by a rotation of angle $\nicefrac{i\pi}{n}$. Note that together $T_0, \dots, T_{n-1}$ cover all chords.
	A flip sequence $F_i$ from $T_i$ to $T_{i+1}$ of length $n$ can be constructed as before, we only have to take special care of the four tangents that replace the long diagonal. Here we use two steps in instead of one. \cref{fig:typeDProof} illustrates an example for $n=4$; the three central triangulations illustrate the modified flips.
	The concatenation of $F_0,\dots, F_{n-1}$ yields a \fh cycle of length $n^2$.
\end{proof}

\section{Facet-Hamiltonian cycles and cluster algebras}
\label{sec:cluster}

In this section, we revisit the constructions of the previous section in the framework of cluster algebras. We give a unified proof of Theorem~\ref{thm:fhgenassoc} relying on standard cluster algebraic tools.

\subsection{Cluster algebras}

Cluster algebras have been defined by Fomin and Zelevinsky in a series of foundational papers~\cite{MR1887642,MR2004457,MR2110627,MR2295199}. For gentle introductions to cluster algebras and further references, we refer the reader to Felikson~\cite{felikson2023ptolemy} and Williams~\cite{MR3119820}. Connections to the Conway-Coxeter friezes are described for instance in the survey paper from Morier-Genoud~\cite{MR3431573}. An extensive treatment of the relation between cluster algebras and triangulated surfaces, beyond the finite cases tackled here, is given by Fomin, Shapiro, and Thurston~\cite{MR2448067}, and Fomin and Thurston~\cite{MR3852257}.
For background on root systems and generalized associahedra, we refer to Bj\"{o}rner and Brenti~\cite{MR2133266}, and Fomin and Reading~\cite{MR2383126}.

Variables of a cluster algebra are grouped into \defn{clusters}, consisting of a pair $(\mathbf{x}, B)$, where $\mathbf{x} = x_1,x_2,\ldots ,x_n$
are called the \defn{cluster variables} and $B=(b_{ij})$ is an $n\times n$ integer matrix called the \defn{exchange matrix}. The matrix $B$ is usually \defn{skew-symmetrizable}.
The initial cluster is called the \defn{seed}. A new cluster can be obtained from any cluster via a \defn{mutation} in the direction $k\in [n]$. The effect
of a mutation in direction $k$ on the cluster variables is the change of the single variable $x_k$ into $x'_k$, satisfying the following relation:
\[
x_kx'_k = \prod_i x_i^{[b_{ik}]_+} + \prod_i x_i^{[-b_{ik}]_+},
\]
where we use the notation $[b]_+ = \max \{b, 0\}$.
The mutation also affects the exchange matrix $B$, which then becomes $B'=(b'_{ij})$ satisfying
\[
b'_{ij} =
\begin{cases}
  -b_{ij} & \text{ if } i = k \text{ or } j = k,\\
  b_{ij} + [b_{ik}]_+ [b_{kj}]_+ - [-b_{ik}]_+[-b_{kj}]_+ & \text{ otherwise.} 
\end{cases}
\]
The cluster algebra is the subring of $\mathbb{Q}(x_1,x_2,\ldots ,x_n)$ generated by the cluster variables of all clusters obtained from the seed by a sequence of mutations.

When the matrix $B$ is skew-symmetric, it can conveniently be represented as a \defn{quiver}, defined as an $n$-vertex directed multigraph, without loops, and without directed 2-cycles. The arcs of a quiver are called \defn{arrows}.
An entry $b_{ij}>0$ is then simply the number of arrows from $i$ to $j$ in the quiver.
The interpretation of a variable mutation in $k\in [n]$ is then as follows:
\[
x_kx'_k = \prod_{\text{arrows from } k \text { to } i} x_i + \prod_{\text{arrows from } i \text { to } k} x_i,
\]
while the quiver mutates according to the following steps:
\begin{itemize}
\item for each subquiver $i\to k\to j$, add an arrow $i\to j$,
\item reverse all arrows incident to $k$,
\item remove all arrows in a maximal set of pairwise disjoint directed 2-cycles.
\end{itemize}

\begin{figure}[htb]
    \centering
    \includegraphics[width=.7\textwidth]{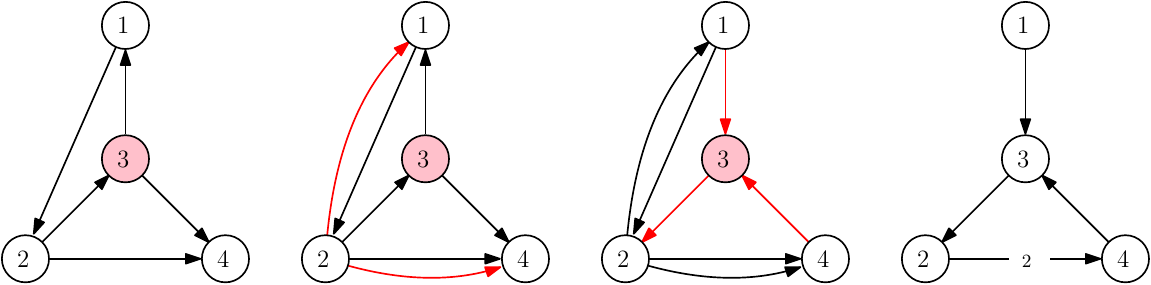}
    \caption{The three steps of mutation of a quiver. Here the quiver on the left is transformed into the quiver on the right by a mutation in 3. Numbers on the arrows indicate multiplicity.}
    \label{fig:quiver}
\end{figure}

One can check that these are indeed the same definitions as using an exchange matrix $B$.
An illustration of the three steps of the mutation of a quiver is given in Figure~\ref{fig:quiver}.

Cluster variables are always rational functions of the seed variables.
An important property of cluster algebra is the \defn{Laurent phenomenon}, that states that the cluster variables are actually always Laurent polynomials in the seed variables.

A cluster algebra is said to be of \defn{finite type} if it has only finitely many seeds, hence the number of distinct cluster variables that can be generated is bounded. A major result of Fomin and Zelevinsky is a classification of the cluster algebra of finite types, matching the Cartan-Killing classification of root systems. In the following statement, the directed graph $\Gamma(B)$ of the exchange matrix $B$ has an arc from $i$ to $j$ if and only if $b_{ij}>0$. 
\defn{Dynkin diagrams} are graphs arising in the classification of Lie algebras and reflection groups.

\begin{theorem}[Fomin-Zelevinsky~\cite{MR2004457}]
  \label{thm:cafinite}
  A cluster algebra is of finite type if and only if it has a seed $(\mathbf{x}, B)$ such that $\Gamma(B)$ is an orientation of a finite type Dynkin diagram.
\end{theorem}

The finite type Dynkin diagrams are well-known and correspond to the four infinite families and the five exceptional irreducible root systems.
We therefore obtain cluster algebras for all finite types $A_n, B_n, C_n, D_n$ and the exceptional types $E_6, E_7, E_8, F_4, G_2$.
The main motivation for introducing cluster algebras of finite type is that they provide a definition of \defn{generalized associahedra} for all finite types.
In cluster algebra of finite types, the \defn{cluster complex} is a simplicial complex whose vertices are cluster variables and maximal simplices are clusters.

\begin{theorem}[Fomin-Zelevinsky~\cite{MR2004457}]
  \label{thm:cagenassoc}
  The cluster complex of a cluster algebra of finite type is the dual of the generalized associahedron of that type.
  In particular, the cluster complexes of type $A$ cluster algebras are duals of associahedra, and those of type $B$ are duals of cyclohedra. 
\end{theorem}

\subsection{Bipartite belts}

We consider the notion of \defn{bipartite belt} defined by Fomin and Zelevinsky~\cite{MR2295199}.
(Note that here we will leave out the so-called \defn{coefficient dynamics}, and only consider the cluster variables.)

An exchange matrix $B$ is said to be \defn{bipartite} if
there exists a function $\varepsilon : [n]\to \{+1, -1\}$ such that
if $b_{ij}$ is positive, then $\varepsilon (i)=+1$ and $\varepsilon (j)=-1$.
If the exchange matrix corresponds to a quiver (hence if it is skew-symmetric), it means that the quiver is bipartite
in the graph-theoretic sense: all arrows are from a positive to a negative vertex.
We will say that a seed $(\mathbf{x}, B)$ is bipartite if $B$ is bipartite.

A simple example of a bipartite quiver is the Dynkin diagram of type $A$, a path \dynkin[Coxeter]A{}, oriented so that sinks and sources alternate. For instance the  matrix
\[
B=
  \left( {\begin{array}{cccc}
  0 & -1 & 0 & 0 \\
  1 & 0 & 1 & 0 \\  
  0 & -1 & 0 & -1 \\  
  0 & 0 & 1 & 0 
  \end{array} } \right)
\]
is the bipartite exchange matrix corresponding to a bipartite orientation of the path on four vertices.

The bipartite belt can be defined on all algebras with a bipartite seed $(\mathbf{x}, B)$.
Note that from the characterization of finite type cluster algebras in Theorem~\ref{thm:cafinite}, and since Dynkin diagrams are bipartite, it holds for all finite type cluster algebras.
Let us denote by $\mu_k$ the mutation in direction $k$.
We define the following two operations:
\[
\mu_+ = \prod_{k: \varepsilon(k)=+1} \mu_k,\ \ \ \  \mu_- = \prod_{k: \varepsilon(k)=-1} \mu_k, 
\]
corresponding respectively to composing the mutations on all vertices on each side of the bipartition.
Given an initial bipartite seed $S_0=(\mathbf{x}_0, B)$, its \defn{bipartite belt} consists of the clusters
\[
S_m = (\mathbf{x}_m, (-1)^m B), m\in\mathbb{Z},
\]
where for $t>0$, $S_t$ is obtained by applying $t$ times one of the operations $\mu_+$ and $\mu_-$, alternately, starting with $\mu_-$:
\[
S_t = \underbrace{\ldots \mu_-\mu_+\mu_-}_{t\text{ factors}}(S_0),
\]
and $S_{-t}$ is defined similarly, but starting with $\mu_+$.
Since for $t>0$ the first mutation is $\mu_-$, we have $x_{1,j} = x_{0,j}$ for all $j$ such that $\varepsilon (j)=+1$,
and in general we have $x_{m+1,j} = x_{m,j}$  for all $j$ such that $\varepsilon (j)=(-1)^m$. 
For $j=(-1)^{m-1}$, on the other hand, we have by definition
\begin{equation}
  \label{eq:bipupdate}
  x_{m-1, j} x_{m+1, j} = 1 + \prod_{i: \varepsilon(i)=-\varepsilon(j)} x_{m,i}^{|b_{ij}|}.
\end{equation}

The following theorem is one of the main result in the fourth paper by Fomin and Zelevinsky~\cite{MR2295199}.

\begin{theorem}[Fomin and Zelevinsky~\cite{MR2295199}]
  \label{thm:bipbelt}
  Consider the bipartite belt constructed from a bipartite exchange matrix $B$ whose graph $\Gamma(B)$ is connected.
  \begin{enumerate}
  \item If the cluster algebra is of finite type, then:
    \begin{enumerate}
    \item the bipartite belt is periodic, and $S_m = S_{m+2(h+2)}$ for all $m\in\mathbb{Z}$, where $h$ is the \defn{Coxeter number},
    \item every cluster variable belongs to a cluster of the bipartite belt,
      \item the denominator vectors establish a bijection between cluster variables and the almost positive roots $\Phi_{\geq -1}$ of the corresponding root system $\Phi$.
    \end{enumerate}
  \item Otherwise, all the elements $x_{m,i}$ in the bipartite belt are distinct.
  \end{enumerate}  
\end{theorem}

\subsection{Bipartite belts in type \texorpdfstring{$A$} and Conway-Coxeter friezes.}

It is convenient to represent the bipartite belt as a \defn{frieze}: For each $m\in\mathbb{Z}$, we represent the cluster variables $x_{m,j}$ such that $\varepsilon (j)=(-1)^m$ in a single column.

In type $A$, this is related to a well-studied family of objects known as the \defn{Conway-Coxeter friezes}~\cite{MR0461269,MR0461270}.
This relation was first described by Caldero and Chapoton~\cite{MR2250855}.
Let the initial bipartite quiver be a path (the type $A$ Dynkin diagram), oriented such that sinks and sources alternate along the path. 
Then Equation~\eqref{eq:bipupdate} governing the update of the cluster variables becomes the Ptolemy relation
\[
ad - bc = 1
\]
for every diamond pattern of variables $a\genfrac{}{}{0pt}{}{b}{c}d$ in the frieze,
where we assume that the top and bottom row all consist of 1s.
This is the rule defining the Conway-Coxeter friezes.
By replacing the initial cluster variables by the number 1 and applying the update rule, we obtain such a frieze.
The bipartite belt in type $A_4$ is represented in Table~\ref{tab:exA4}, and a frieze obtained by evaluating the cluster variables with $x_i=1$ for all $i$ is given in Table~\ref{tab:exfrieze}.

\begin{table}[t]
\caption{\label{tab:exA4}The bipartite belt in type $A_4$. For every 4-tuple of variables forming a diamond $a\genfrac{}{}{0pt}{}{b}{c}d$, and assuming top and bottom rows consisting of 1s, the relation $ad - bc = 1$ holds.}
\noindent\makebox[\textwidth]{%
  $
   \begin{array}{ccccccccc}
  m= -1   & 0 & 1 & 2 & 3 & 4 & 5 & 6 & 7 \\
  \hline \\
   x_1 &  & \frac{1+x_2}{x_1} & & \frac{x_1 x_3 + x_2 x_4 + 1}{x_2 x_3} & & \frac{x_3 + 1}{x_4} & & x_4 \\
 & x_2 &  & \frac{x_2^2 x_4 + x_1 x_3 + x_2 x_4 + x_2 + 1}{x_1 x_2 x_3} & & \frac{x_1 x_3^2 + x_1 x_3 + x_2 x_4 + x_3 + 1}{x_2 x_3 x_4}  & & x_3 &  \\
   x_3 &  & \frac{x_2x_4 + 1}{x_3} & & \frac{x_1 x_3^2 + x_2^2 x_4 + x_1 x_3 + x_2 x_3 + x_2 x_4 + x_2 + x_3 + 1}{x_1 x_2 x_3 x_4}  &   & \frac{x_1 x_3 + 1}{x_2} & & x_2 \\
   & x_4 &  & \frac{x_2 x_4 + x_3 + 1}{x_3 x_4} & & \frac{x_1 x_3 + x_2 + 1}{x_1 x_2} & & x_1 & \\
   \\
   \hline
   \end{array}
   $}
\end{table}

\begin{table}[t]
\centering
\caption{\label{tab:exfrieze}A frieze obtained from the type $A_4$ bipartite belt.}
    $
   \begin{array}{ccccccccc}
       \hline
   1 &  & 2 & & 3 & & 2 & & 1 \\
 & 1 &  & 5 & & 5 & & 1 &  \\
   1 &  & 2 & & 8 & & 2 & & 1 \\
   & 1 &  & 3 & & 3 & & 1 & \\
   \hline
   \end{array}
   $
\end{table}

\begin{figure}[htb]
\centering
  \includegraphics[scale=.4]{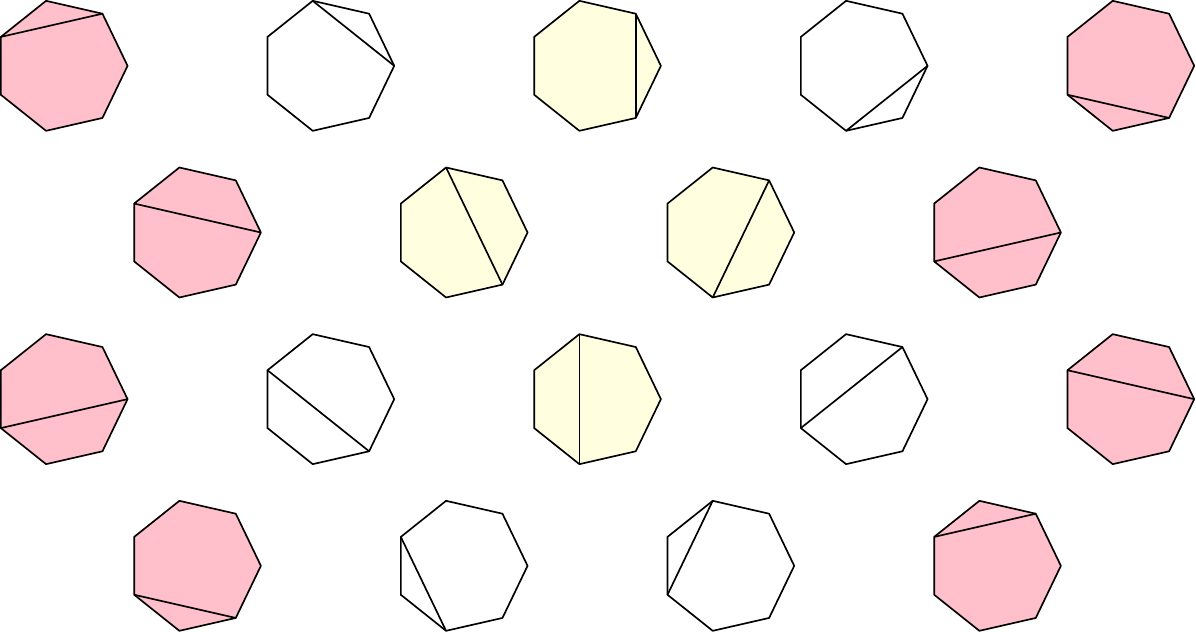}
  \caption{\label{fig:frieze}The type $A_4$ frieze with diagonals of a heptagon.}
  \end{figure}

Similarly, we can replace every cluster variable by a diagonal of the $(n+3)$-gon, such that the initial seed maps to a zigzag triangulation of the polygon. The arrows in the seed quiver can be interpreted as pairs of diagonals forming a clockwise angle.
An application of either $\mu_-$ or $\mu_+$ then corresponds to a rotation of angle $2\pi/(n+3)$ of every other diagonal, and the successive applications of both yield a rotation of the whole triangulation. An illustration is given on \cref{fig:frieze}. 

The update of the cluster variables for the diamond patterns $a\genfrac{}{}{0pt}{}{b}{c}d$ then corresponds to a diagonal flip in the triangulation, and the relation resembles the Ptolemy relation on the length of diagonals of a quadrilateral. (This interpretation is in fact correct in a hyperbolic setting.) By setting to 1 the variables corresponding to another triangulation, we obtain another frieze, the second row of which is the \defn{quiddity} of the triangulation: the number of triangles incident to each successive vertex of the polygon.
We refer to the original papers from Conway and Coxeter~\cite{MR0461269,MR0461270}, and to Morier-Genoud’s survey~\cite{MR3431573} for details on these beautiful structures.

\subsection{Facet-Hamiltonian cycles from friezes}
\label{sec:friezestofh}

We now explain how we can easily extract many \fh cycles on generalized associahedra from the bipartite belts of the corresponding cluster algebras of finite types.

Given a bipartite belt in finite type, we construct a \fh cycle by starting with the initial seed, and performing the mutations $\mu_k$ for all $k$ such that $\varepsilon(k)=-1$, the composition of which yields $\mu_-$, then the mutations $\mu_k$ for all $k$ such that $\varepsilon(k)=+1$, the composition of which yields $\mu_+$.
This is mutating the initial seed $S_0$ into $S_2$, and visiting all cluster variables $\mathbf{x}_2$. Iterating, we will obtain all cluster variables, which, from Theorem~\ref{thm:cagenassoc}, correspond to facets of the corresponding generalized associahedron.
From Theorem~\ref{thm:bipbelt}, this will yield a cycle that passes through every cluster variable exactly once, hence every facet of the generalized associahedron exactly once, as desired.
This directly proves our main result.

\assoc*

\begin{figure}[htb]
\centering
  \includegraphics[page=3]{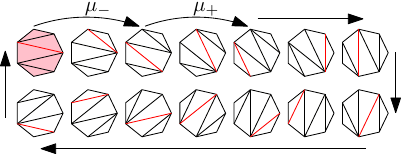}
  \caption{\label{fig:cycleA}A \fh cycle on the type $A_4$ associahedron.}
\end{figure}

\enlargethispage{12pt}
As explained in the previous section, for the three main types $A$, $B/C$, and $D$, the bipartite belts have an interpretation in terms of triangulation models.
An example of cycle obtained from the bipartite belt of type $A_4$ is given in Figure~\ref{fig:cycleA}. One can check that  
it is precisely the cycle described in the second proof of Section~\ref{sec:Aassoc} and illustrated in \cref{fig:associahedron}.

The cluster algebras of type $B$, with Dynkin diagrams  of the form \dynkin B{}, do not have skew-symmetric exchange matrices, hence are not quiver cluster algebras, but there is a natural interpretation of the bipartite belt as rotation of diagonals in symmetric triangulations of a $2n$-gon. Again, the second cycle described in Section~\ref{sec:Bassoc} can be read from the type $B$ bipartite belt.

\begin{figure}[htb]
\centering
  \includegraphics[scale=.4, page=2]{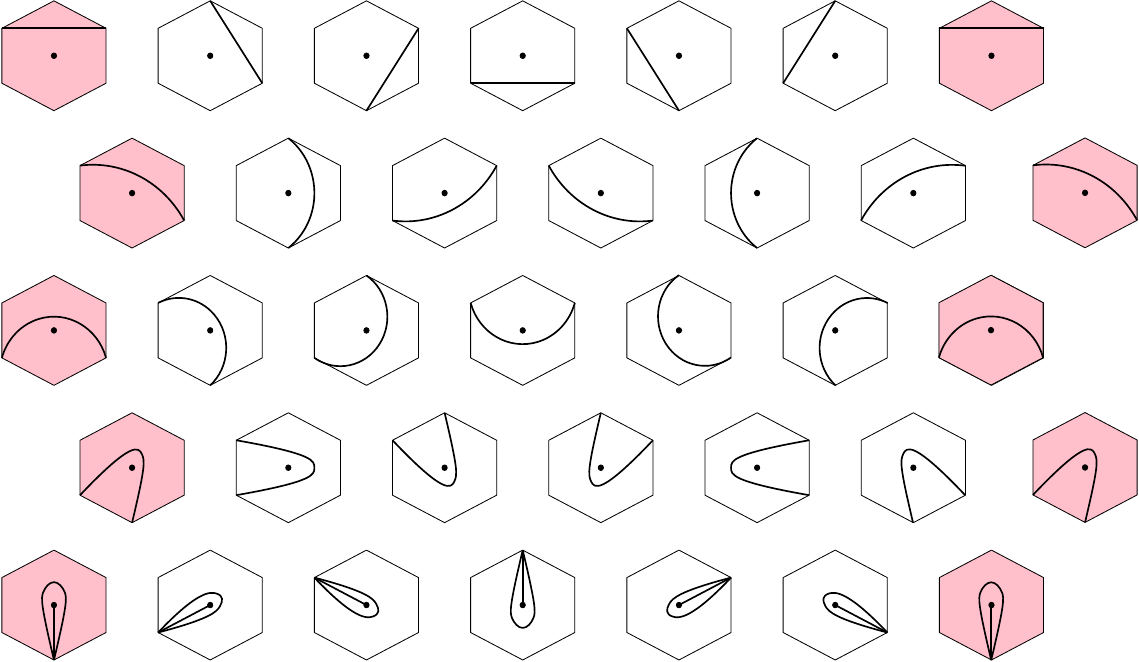}
  \caption{\label{fig:triangD}The type $D_6$ frieze with arcs of a punctured hexagon.}
\end{figure}

Bipartite exchange matrices for the type $D$ have diagrams of the form \dynkin D{}.
Orienting the arrows such that sinks and sources alternate yield bipartite quivers. In Figure~\ref{fig:triangD}, we show the obtained frieze using arcs in a triangulation of a punctured $n$-gon as models for the cluster variables.
This triangulation model is described by Fomin, Shapiro, and Thurston~\cite{MR2448067}, and is obtained from the symmetric pseudo-triangulation model by Ceballos and Pilaud~\cite{CP16} described in the previous section by folding the triangulation around its center of symmetry.

In fact, bipartite belts directly yield rhombic strips in a suitable poset. All the \fh cycles on generalized associahedra of types $A$, $B/C$, and $D$ of \cref{sec:assoc} can be read from the corresponding friezes, by considering initial clusters that do not correspond to bipartite quivers, but are constructed from the frieze as a path from the top to the bottom row. In the types $A$, $B/C$ and $D$, these cluster correspond to triangulations in which all triangles have at least one edge which is an edge of the polygon. The obtained cycles include for instance those described on Figures~\ref{fig:associahedron1} and~\ref{fig:cycle3}.

Bipartite belts are used in software packages for cluster algebras, with the purpose of generating all cluster variables, even in cluster algebras of infinite types~\cite{MR2862990}. Since from Theorem~\ref{thm:bipbelt}, the cluster variables of the bipartite belt are all different, it provides a simple way to generate arbitrarily many in a systematic fashion. However, it is not true in general that all cluster variables lie in the bipartite belt, so some cluster variables may never be generated this way. Still, the idea has further applications, see for instance Assem, Reutenauer, and Smith \cite{MR2729004}, and Pallister~\cite{MR4496491} for applications to affine quivers. 

\section{Facet-Hamiltonian paths and cycles in graph associahedra}
\label{sec:grassoc}

We first recall the definition and properties of graph associahedra in Section~\ref{sec:defs}. In Section~\ref{sec:tools}, we introduce a number of simple tools and definitions that will be helpful in the remainder. \cref{sec:csgfw,sec:caterpillar,sec:completeBipatite} give the proofs of \cref{thm:grassoc} and \cref{thm:bipcat}.

\subsection{Graph associahedra}
\label{sec:defs}

For this entire section, let $G=(V,E)$ be a simple connected graph with $n:=|V|$.

\begin{figure}[htb]
    \centering
    \begin{subfigure}[t]{0.15\textwidth}
    \centering
        \includegraphics[scale = 0.6]{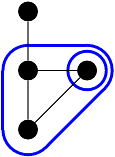}
        \subcaption{nested}
        \label{fig:compatible1}
    \end{subfigure}
    \hfill
    \begin{subfigure}[t]{0.25\textwidth}
    \centering
        \includegraphics[scale = 0.6]{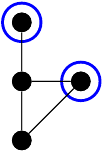}
        \subcaption{disjoint, non-adjacent}
        \label{fig:compatible2}
    \end{subfigure}
    \hfill
    \begin{subfigure}[t]{0.25\textwidth}
    \centering
        \includegraphics[scale = 0.6]{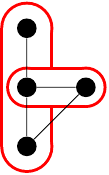}
        \subcaption{neither disjoint nor nested}
        \label{fig:notcompatible1}
    \end{subfigure}
    \hfill
    \begin{subfigure}[t]{0.25\textwidth}
    \centering
        \includegraphics[scale = 0.6]{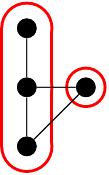}
        \subcaption{disjoint but adjacent}
        \label{fig:notcompatible2}
    \end{subfigure}
    \caption{Compatible pairs of tubes in ({\sc a}) and ({\sc b}), and non-compatible pairs of tubes in ({\sc c}) and ({\sc d}).}
    \label{fig:compatible}
\end{figure}

We recall the definition of tubes and tubings given in \cref{sec:intro}. 
A \defn{tube} of $G$ is a nonempty proper subset $t\subset V$ such that the induced subgraph $G[t]$ is connected. 
Two tubes $t_1$ and $t_2$ are \defn{compatible} if one of the following conditions is fulfilled:
\begin{itemize}
    \item They are \defn{nested}: either $t_1\subset t_2$ or $t_2\subset t_1$.
    \item They are \defn{non-adjacent}: $G[t_1\cup t_2]$  is not connected.
\end{itemize}
(Note that non-adjacent tubes are necessarily disjoint.) 
\cref{fig:compatible} depicts examples of pairs of tubes.

A \defn{tubing} is a collection of pairwise compatible tubes. It is called \defn{nested} if all pairs of its tubes are nested.

\begin{figure}[htb]
    \centering
    \includegraphics[scale = 0.6]{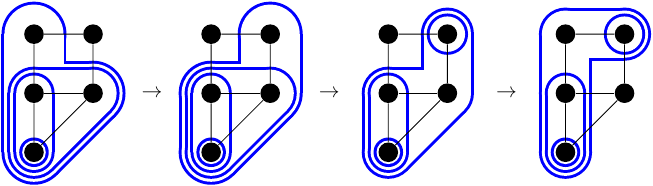}
    \caption{An example of a  flip sequence. The first two tubings are nested.}
    \label{fig:flips}
\end{figure}

The graph associahedron $\mathcal A(G)$ is the polytope whose face lattice is obtained by ordering the tubings of $G$ in reverse inclusion order. The vertices of $\mathcal A(G)$ therefore correspond to the maximal tubings of $G$, and the facets correspond to tubings of size one, hence to single tubes.  This polytope is simple because every maximal tubing of a connected $n$ vertex graph consists of $n-1$ tubes and each tube of a maximal tubing can be replaced by a unique distinct tube. We say that this tube is being \defn{flipped}, see  \cref{fig:flips} for examples. Edges of $\mathcal A(G)$ are in bijection with pairs of maximal tubings that differ by a single flip.

In this paper, we  only consider graph associahedra of connected graphs. This is because the associahedron $\mathcal A(G)$ of a graph $G$ with multiple connected components $C_1,\dots, C_k$ is the cartesian product of the graph associahedra of the components: $\mathcal A(G)=\mathcal A(C_1)\times\dots \times\mathcal A(C_k)$. Thus, $\mathcal A(G)$ has a  \fh cycle (or path) if and only if each  $\mathcal A(C_i)$ has a \fh cycle (or path). 

\subsection{Graph associahedra and triangulations.}

Graph associahedra generalize many well-known families of polytopes, including some that were already discussed in previous sections. In particular, (type A) associahedra are graph associahedra in which the graph is a path. 

\begin{figure}[htb]
	\centering
	\begin{subfigure}{.45\textwidth}
		\centering
		\includegraphics[scale = 0.5]{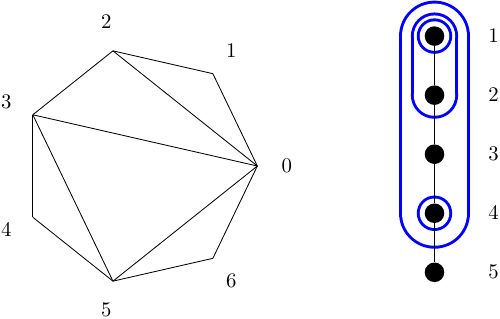}
		\subcaption{}
		\label{fig:bijectionA}
	\end{subfigure}\hfil
	\begin{subfigure}{.4\textwidth}
		\centering
		\includegraphics[scale = 0.5]{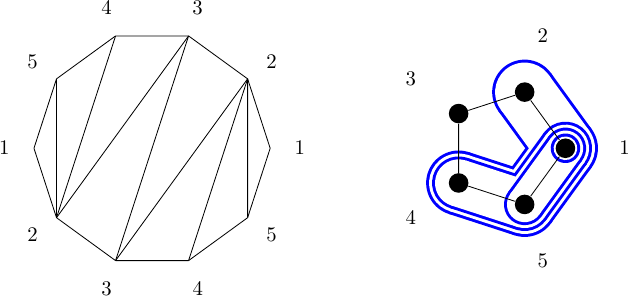}
		\subcaption{}
		\label{fig:bijectionB}
	\end{subfigure}
	\caption{Bijections between ({\sc a}) diagonals of the $(n+2)$-gon and tubes of $P_n$ and ({\sc b}) diagonals of the $2n$-gon and tubes of $C_n$.}
	\label{fig:bijection}
\end{figure}

The bijection between triangulations of a convex $(n+2)$-gon and maximal tubings of a path is obtained as follows. Label the vertices of the $(n+2)$-gon from 0 to $n+1$ in, say, clockwise order, and the vertices of the path from 1 to $n$. For an example consider \cref{fig:bijectionA}. Now consider a triangulation, and for a diagonal of the triangulation connecting the vertices $i$ and $j$, with $j>i+1$, include the tube $\{i+1,\ldots ,j-1\}$. It is easy to see that this collection of tubes is a maximal tubing, and that the map is bijective. The map generalizes to an order-preserving map between the non-crossing sets of diagonals and (not necessarily maximal) tubings of the path.

Similarly, there is a bijection between the tubes (respectively, maximal tubings) of the $n$-cycle~$C_n$  and the diagonals (respectively, symmetric triangulations) of the $2n$-gon, implying that type B associahedra are associahedra of cycles. Label the corners of the $2n$-gon by $1,\dots, n, 1,\dots, n$ as in \cref{fig:bijectionB}. Here a $k$-tube of $C_n$ is mapped to the diagonals that span the corresponding $k$ elements: it connects the predecessor and the successor. For example the tube $\{1,2,3\}$ corresponds to the diagonals $\{n,4\}$. It is easy to see that all bistars correspond to permutations that are cyclic shifts of $1,\dots, n$.
In the case where the graph is complete, every maximal tubing consists of pairwise nested tubes, and every flip preserves the size of the flipped tube. Flipping the tube of size $k$ is equivalent to swapping the elements in positions $k$ and $k+1$ in the corresponding permutation, which is an adjacent transposition. In this case the graph associahedron is the permutahedron.

Finally, note that the stellohedra, the graph associahedra of stars, encode regular triangulations of a family of point sets known as the \defn{mother of all examples}~\cite{LRS10}. We refer the reader to Ceballos and Pilaud~\cite{CP16} for more details.

\subsection{Facet-Hamiltonian paths/cycles in graph associahedra.}

The vertices incident to a facet of the graph associahedron $\mathcal A(G)$ are all the maximal tubings containing a certain tube. Thus, for a \fh cycle or path, we seek sequences of maximal tubings of $G$ such that every tube appears in an interval along this sequence and every adjacent pair in the sequence is related by a flip. See \cref{fig:flips} for an example of a part of such a sequence. One can ask for a \fh cycle or path to have the additional property that all its tubings are nested, which will be fulfilled for many of the constructions which follow. In this case the cycle or path can be represented as a rhombic strip similar to the constructions seen in Figures \ref{fig:rhombic3}, \ref{fig:rhombic4} and \ref{fig:rhombic3A}. In this case we seek a rhombic strip which is a subgraph of the Boolean lattice obtained by deleting all subsets which do not correspond to tubes. This is the inclusion poset on the tubes. See \Cref{fig:cycle rhombic representation} for an example.

\subsection{Tools}
\label{sec:tools}

We now introduce a few operations that will be essential in our constructions of \fh cycles and paths.

In what follows, we will refer to nested maximal tubings simply as \defn{nested tubings}. A nested tubing has a unique tube consisting of only one vertex, which we call the \defn{\kernel} of the tubing. A nested tubing corresponds to a permutation by ordering the vertices according to the number of tubes they are contained in, in decreasing order, so that the \kernel is the first element of the permutation. 

\begin{figure}[htb]
    \centering
    \includegraphics[scale = 0.6]{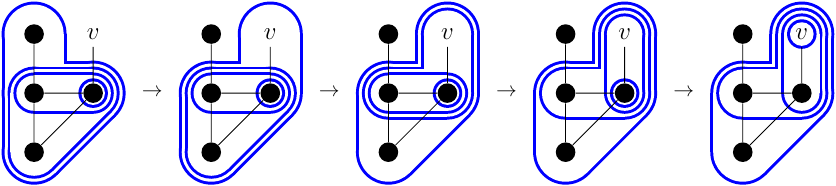}
    \label{fig:absorbing}
    \caption{Vertex $v$ is being absorbed into the tubing. In this example this induces four flips.}
\end{figure}

Let us consider a graph $G^{+v}=(V\cup\{v\},E')$ with $E\subset E'$, obtained by adding $v$ to $G$, together with some new edges incident to $v$. 
Let $T$ be a nested tubing of $G$ such that its \kernel is adjacent to $v$. There is a labeling of the vertices such that $T$ is associated with the permutation $\pi = 1,\ldots,n$, with \kernel~1. We define the procedure \defn{absorbing} $v$ into $T$ as follows. 
Let us define a new nested tubing $T' \coloneqq T\cup\{V\}$ associated with the permutation $\pi_{n+1} = 1,\ldots,n,v$. We now observe that flipping the tubes of $T'$ in descending order, thus iteratively replacing the tube $[k]$ by $[k-1]\cup\{v\}$, gives a new nested tubing associated with $\pi_k$ in each step. Thus, in terms of the associated permutations, these flips correspond to the adjacent transpositions in descending order. The final nested tubing has \kernel $v$ and the permutation $\pi_1 = v,1,\ldots,n$. \cref{fig:absorbing} illustrates this process.
The process of absorbing $v$ can be reversed by flipping all tubes in ascending order. We will refer to this procedure as \defn{expelling} $v$.

\begin{restatable}{lemma}{absorbtion}\label{lem:absorption}
    Let $T$ be a nested tubing of $G$ with permutation $\pi$ and \kernel $c$, and $G^{+v}$ a supergraph of $G$ obtained by adding a vertex $v$ such that $cv\in E(G^{+v})$. Then $T'=T\cup\{V\}$ is a nested tubing and the process of absorbing $v$ produces a valid path in the associahedron~$\mathcal A(G^{+v})$ consisting of the nested tubings associated to $\pi_{n+1},\dots, \pi_1$. Furthermore, each tube introduced in this path is present in the last tubing associated with the permutation $\pi_1$.
\end{restatable}

 \begin{proof}
    $T'$ is indeed a nested tubing because every tube in $T$ is contained in $V$. Given a labeling of $G$ such that $T$ is associated with permutation $\pi = 1,\ldots,n$, in step $k$ we replace tube the $[n-k+1]$ with $[n-k]\cup\{v\}$. This new tube is compatible with all other tubes because it is either a subset or a superset in each case. Furthermore, $[n-k]\cup\{v\}$ is actually a tube because $v$ is connected to the \kernel $c=1$, hence $[n-k]\cup\{v\}$ is connected in $G^{+v}$. The tubes present in the final tubings are the prefixes of the permutation $\pi_1$.
\end{proof}

Now let us consider a path $P$ in $\mathcal A(G)$ such that for each maximal tubing $T$ in $P$ the following properties hold: (i) $T$ is nested, (ii) the \kernel of $T$ is adjacent to $v$ in the supergraph~$G^{+v}$. Let $P^{+v}$  be the path obtained as follows: For each maximal tubing $T$ of $P$, add $v$ to every tube of $T$, and add the tube $\{v\}$ to $T$.
In other words, a tubing in $P$ with  permutation $\pi$ corresponds to a tubing in $P^{+v}$ with permutation $\pi_1$.

\begin{restatable}{lemma}{pathabsorbtion}\label{lem:P0}
    $P^{+v}$ is a well-defined path in $\mathcal A(G^{+v})$.
\end{restatable}

\begin{proof}
     We have seen in \cref{lem:absorption} that every tubing in $P^{+v}$ is well-defined. Let $T_1$ and $T_2$ be consecutive tubings in $P$. Then there are tubes $t_1$ and $t_2$ such that $T_2 = T_1\backslash\{t_1\}\cup\{t_2\}$. Let $T_1'$ and $T_2'$ be the corresponding adjacent tubings in $P^{+v}$. Consider $t_1'=t_1\cup\{v\}$ and $t_2'=t_2\cup\{v\}$. Then $T_2'=T_1'\backslash\{t_1'\}\cup\{t_2'\}$. This shows that $T_1'$ and $T_2'$ are connected by an edge in $\mathcal A(G^{+v})$.
\end{proof}

Combining the above two lemmas and ideas from the proof of \cref{thm:perm}, we can construct \fh cycles in graphs obtained by adding a \defn{universal vertex}; we call a vertex \defn{universal}  if it is adjacent to every other vertex of a given graph. 

\begin{lemma}\label{lem:UniversalVertex}
	Let $P$ be a \fh path of the graph associahedron $\mathcal{A}(G)$. 
        Let $G^+$ be the graph obtained from $G$ by adding a universal vertex $v$ to $G$. 
        Let $V=[n]$ be a labeling of the vertices.
	\begin{enumerate}
		\item If $P$ ends in the nested tubing associated with $\pi = 1,2,\dots,n$, then $\mathcal A(G^+)$ has a \fh path ending in the nested tubing associated with $v,1,\ldots,n$.
		\item If $P$ starts in $n,1, \ldots, n-1$ and ends in $1, 2, \ldots, n$, then $\mathcal A(G^+)$ has a \fh cycle.
	\end{enumerate}
\end{lemma}
\begin{proof}
    We follow the basic idea of the proof of \cref{thm:perm}. Let $T$ be the final nested tubing of~$P$, and let $Q$ be the sequence of nested tubings obtained by absorbing $v$ into $T$. From \cref{lem:absorption}, $Q$ is a path in~$\mathcal A(G^+)$.
    Let $P'$ be the path found in the proof of \cref{thm:perm}.
    From Lemma~\ref{lem:P0}, the sequence~$P'^{+v}$ is also a path in $\mathcal A(G^+)$.
    Then we claim that under the assumption in (i), the concatenated sequence~$PQ(P'^{+v})$ is a \fh path in $\mathcal A(G^+)$. In $P$ we see all tubes not containing~$v$, 
    and in $P'^{+v}$ we see all tubes containing~$v$. In $Q$ we only encounter those tubes which are either already present in the end of $P$ or in the beginning of $P'^{+v}$. 

    Now suppose that (ii) is fulfilled. Let $Q\inv$ be the path obtained by expelling $v$ from the tubing at the end of $P'^{+v}$. By construction and Lemma~\ref{lem:absorption} again, $PQ(P'^{+v})Q\inv$ is a \fh cycle.
\end{proof}

\begin{remark}
    If a graph $G$ fulfills condition (i) or (ii)  of \cref{lem:UniversalVertex}, then $G^+$ fulfills the same condition. Thus, by repeated application of \cref{lem:UniversalVertex}, we may add universal cliques.
\end{remark}

The proofs of the following results rely heavily on absorption and expulsion. There we will not always absorb a vertex that is not contained in any tube yet, but start the procedure from the middle. Likewise, vertices will be expelled that are not the kernel of their nested tubing. In the context it will always be clear what is being done.

\subsection{Complete split graph, fan, and wheel associahedra}
\label{sec:csgfw}

    We now present our construction of \fh paths and cycles for complete split graphs, fans and wheels. As a first step, we revisit the \fh cycles for type A and B/C associahedra that we have seen in \cref{sec:assoc}. We give equivalent proofs  in the language of graph associahedra. By doing so, we also find \fh paths  where we pay close attention to the start and end of these paths.

\begin{restatable}{proposition}{goodPathforStarPathCycle}\label{lem:goodPathforStarPathCycle}
	The graph associahedron $\mathcal A(G)$ of a graph $G$ on $n$ vertices has a \fh path that starts in $1, \ldots, n$  and ends in $2,\ldots, n, 1$, if $G$ is
	\begin{enumerate}
		\item the path $P_n$ (where the vertices are labeled $1,2,\ldots,n$ along the path)
		\item the cycle $C_n$ (where the vertices are labeled $1,2,\ldots,n$ along the cycle).
  	    \item \label{item:i} the star $S_n$ (where the center has label 1)\footnote{Here, $S_n$ is the star with $n$ vertices and $n-1$ leaves. We use this non-standard definition is to make the dimensions of all polytopes mentioned in the proposition the same.},
	\end{enumerate}
    Furthermore, in these cases $\mathcal A(G)$ has a \fh cycle.
\end{restatable}

\begin{figure}[htb]
    \centering
    \includegraphics[scale = 0.25]{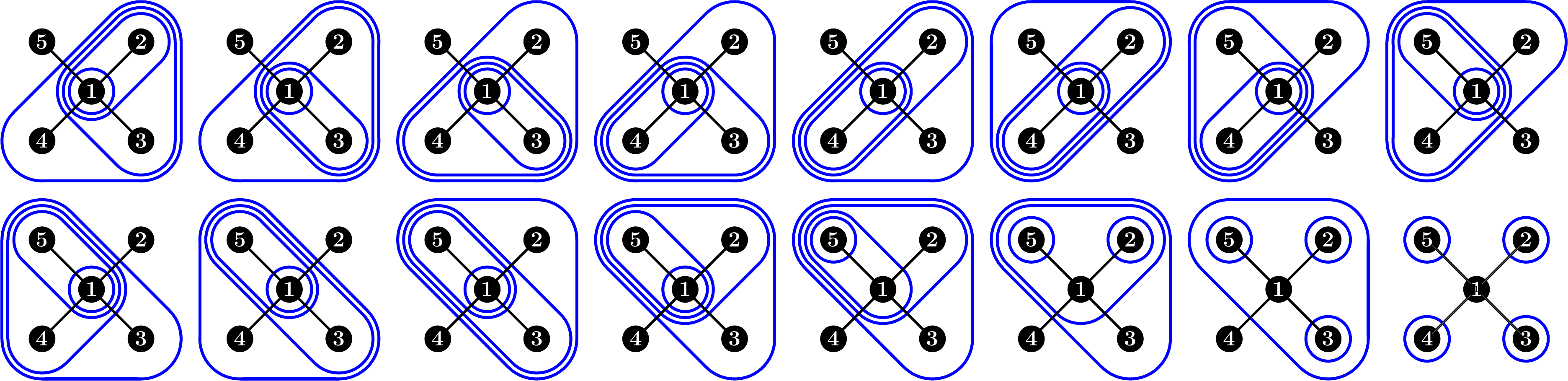}
    \caption{The \fh path obtained on $\mathcal{A}(S_5)$ by the construction in the proof of \cref{lem:goodPathforStarPathCycle}(\ref{item:i}).}    \label{fig:star}
\end{figure}

\begin{proof}
    Let $G=P_n$ and label the vertices along the path by $[n]$. Then the tubes are the intervals of~$[n]$.
    First we give a general construction for a \fh cycle of $\mathcal A(P_n)$. Start with the nested tubing $T=\{t_1,\ldots,t_n\}$ associated with $1,\ldots,n$. Now we apply the following algorithm. In each step we find a tube $t_k$ such that either $k=n-1$ and $t_k = [n-1]$ or the only element in $t_{k+1}\backslash t_k$ is the maximum of $t_{k+1}$. (Note that there may be many candidates; for our depicted cycles we flip the tubes by increasing size.) Such a tube always exists if $T\ne\{\{n\},\ldots,\{2,\ldots,n\}\}$. Flip $t_k$ and update $T$. Repeat this until $T=\{\{n\},\ldots,\{2,\ldots,n\}\}$. Now flipping every tube once more gives a \fh cycle. This can be seen as follows. Draw the vertices of $P_n$ in a vertical line, labeled $1,\ldots,n$ from top to bottom. Each flip described above consists of shifting a tube of $T$ one step downwards. After applying the algorithm each tube has been shifted all the way down and has hence passed over every possible interval of its respective size. Furthermore, once an interval has been left behind, it is clearly not seen again because we never shift upwards. This shows that the algorithm produces a \fh path. The last step then joins the ends of this path as can be seen in \cref{fig:pathCycleAmain}. Now the claimed \fh path can be obtained by breaking this cycle open in the appropriate spot.

    \begin{figure}[htb]
    \centering
        \includegraphics[scale = 0.5]{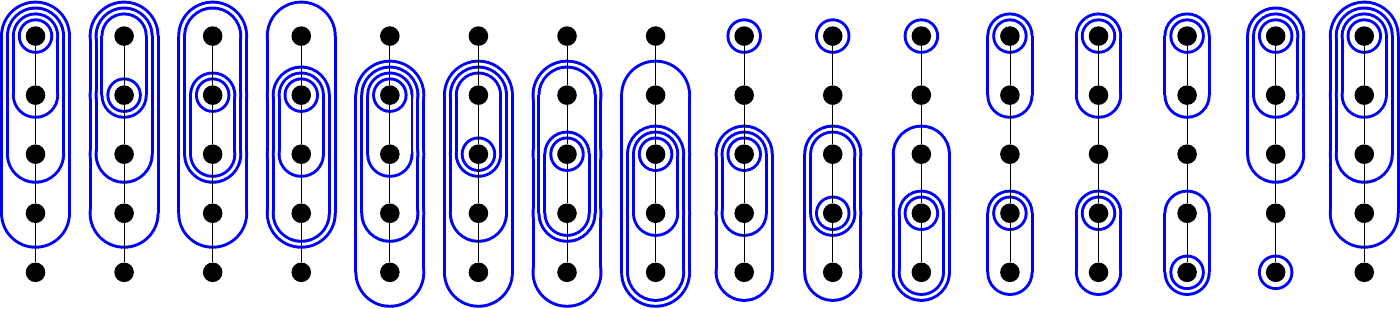}
        \caption{A \fh cycles for the path $P_5$. For $n=4$, this construction yields the cycle 
        from \cref{fig:associahedron3}.}
        \label{fig:pathCycleAmain}
    \end{figure}

    Now, we consider the case where $G=C_n$. Label the vertices along the cycle. Start again in the nested tubing associated with $1,\ldots,n$. We apply a similar idea as for the path. In a phase, we flip the tubes in the order of increasing size each once. On the level of the permutations that means that element 1 is shifted to the right in the first phase. Thus, after the first phase, we obtain the tubing associated with $2,\dots, n, 1$. In general, a phase corresponds to a cyclic left-shift in the permutation, and after $n$ phases we are back where we started. This can be seen as shifting the intervals along the cycle in (for example) a clockwise direction, so the same reasoning as for the path shows that this gives a \fh cycle. Again, breaking the cycle open in the right spot gives the path we are after. Here it suffices to remove the last phase to obtain a path ending in $n,1,\ldots,n-1$. \cref{fig:pathCycleBb} illustrates the result for $C_4$ and \cref{fig:cycle rhombic representation} a \fh cycle of $\mathcal{A}(C_6)$ as a rhombic strip.

       \begin{figure}[htb]
    \centering
        \includegraphics[scale = 0.5]{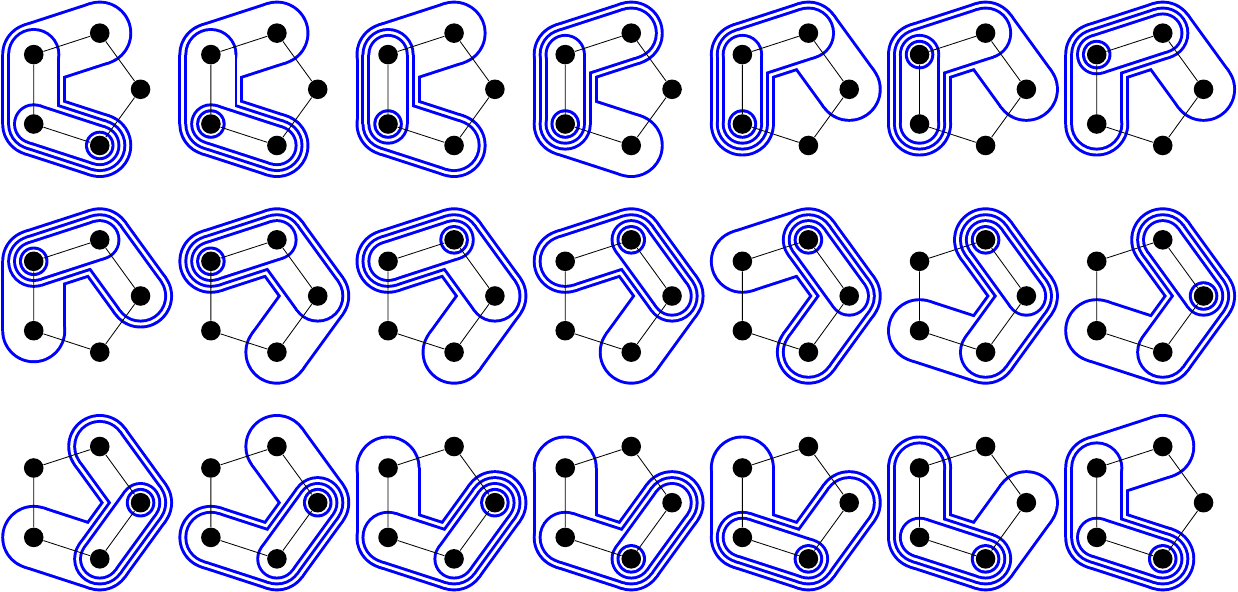}
        \caption{A \fh cycle for the cycle $C_5$. For $n=4$, this construction yields the cycle from \cref{fig:cyclo3}.}
        \label{fig:pathCycleBb}
    \end{figure}

      \begin{figure}[hbt]
        \centering
        \includegraphics[scale = 0.7]{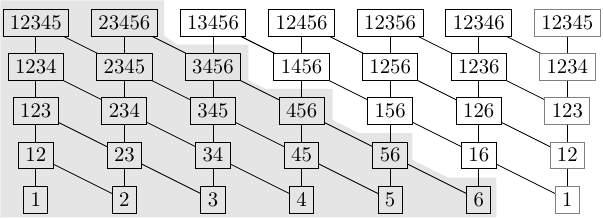}
        \caption{A \fh cycle of $\mathcal{A}(C_6)$ as a rhombic strip. The cycle is assumed to be labeled $1,\ldots,6$ around. The shaded area encodes \fh paths on $\mathcal{A}(P_6)$.}
        \label{fig:cycle rhombic representation}
    \end{figure}

    Finally, let us consider the case where 
    $G=S_n$ with center $1$. Note that the tubes containing $1$ consist of arbitrary subsets of $[n]$ that contain $1$. Let $P'$ be the \fh path of $\mathcal A(K_{n-1})$ from \Cref{lem:fhpaths} where we use the labels $\{2,\ldots,n\}$ on $K_{n-1}$. Let $P$ be the path on $\mathcal{A}(G)$ obtained from $P$ by appending $1$ at the beginning of every permutation of $P'$. Then along $P$ we see all tubes of $G$ that contain the center. Lastly, expelling 1 yields the tubing $\{\{2\}, \ldots, \{n\}\}$ and hence a \fh path. It is clear that we can flip the tubing at the end of this path to any nested tubing, so this path can be closed to a \fh cycle. \cref{fig:star} illustrates the result for $S_4$.
     
     It remains to show that the permutations match our claims. First of all, the beginning and the end of $P$ really are nested permutations. The beginning is associated with $1,\ldots,n$ by construction. $P'$ shifts the elements $2,\ldots,n$ to the right. Furthermore, the end of $P$ consists of a tubing where 1 is not the \kernel. So 1 must be in the second position of the permutation. It follows that the permutation associated with the end of $P$ is $n,1,\ldots,n-1$. Reversing $P$ hence gives the desired path. The ends of this path can be joined by expelling $n$ to obtain a \fh cycle.
\end{proof}

The star, the path and the cycle all fulfill the second condition of \cref{lem:UniversalVertex}. This gives \fh cycles for complete split graphs, fans and wheels. This completes the proof of \cref{thm:grassoc} (together with \cref{thm:perm,lem:asso,lem:cyclo}).
\begin{corollary}
 The graph associahedron $\mathcal A(G)$ of a graph $G$ has a \fh cycle if $G$ is  a fan, a wheel, or a complete split graph.
\end{corollary}
\begin{proof}
Note that a fan and a wheel contain exactly one universal vertex. After its removal, we obtain a path or cycle respectively. A complete split graph contains many universal vertices; removing all but one, yields a star. \cref{lem:goodPathforStarPathCycle} shows that the star, the path and the cycle all fulfill condition (ii) of \cref{lem:UniversalVertex}. This completes the proof.
\end{proof}

\graphAssociahedra*

\subsection{Caterpillar associahedra}\label{sec:caterpillar}

    A \defn{caterpillar} is a tree that becomes a path if all leaves are removed. 
    Let $G$ be a graph and let $P=(T_1,\ldots,T_k)$ be a ordered subset of $\mathcal A(G)$. For a tube $t$ of $G$,  $P_t$ denotes the subset of $P$ such that $t\in T_i\iff i\in P_t$. We say that $P_t$ is \defn{connected}, if these tubings form an interval along $P$ that is not empty.

\begin{lemma}\label{lem:caterpillar}
    Graph associahedra of caterpillars have a \fh path.
\end{lemma}

\begin{proof}
    A spine of $G$ is a path obtained by deleting leaves.
    We denote the vertices of a spine of $G$ by $s_1,\ldots,s_n$ for $n\in\N$ and the leaves of $s_i$  by $\ell_{i,j}$ for $i=1,\ldots,n$ and $j=1,\ldots,m_i$ where $m_i$ is the number of leaves adjacent to $s_i$. Note that we allow $s_n$ to be a leave. For a schematic illustration, consider \cref{fig:catA}.

    \begin{figure}[htb]
	\centering
	\begin{subfigure}{.45\textwidth}
		\centering
		\includegraphics[page=1]{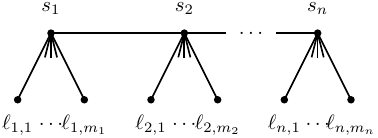}
		\subcaption{}
		\label{fig:catA}
	\end{subfigure}\hfil
	\begin{subfigure}{.4\textwidth}
		\centering
		\includegraphics[page=3]{figures/AssoCaterpillar.pdf}
		\subcaption{}
		\label{fig:catB}
	\end{subfigure}
	\caption{Illustration for the proof of \cref{thm:bipcat} for caterpillars.}
	\label{fig:cat}
\end{figure}

    We apply induction on the number of vertices of $G$. The induction hypothesis is as follows: For any caterpillar $G$ with vertex labeling as above the graph associahedron $\mathcal{A}(G)$ contains a path $P$ with the following properties:
    \begin{itemize}
        \item All tubings of $P$ are maximally nested.
        \item $P$ starts in the maximally nested tubing associated to
        \[
            s_1,\ell_{1,1},\ldots,\ell_{1,m_1},\ldots,s_n,\ell_{n,1},\ldots,\ell_{n,m_n}.
        \]
        \item $P_t$ is connected exactly for the tubes $t$ of $G$ that do not consist of only one leave of $G$.
        \item The path $P$ splits into two paths $P_0$ and $X$ such that all tubes in all tubings of $X$ contain $s_n$ and all such tubes appear along $X$.
    \end{itemize}
    
    For the induction base we consider the caterpillar on two vertices with labels $s_1$ and $\ell_{1,1}$. Then the path $P$ consist of the tubing corresponding to $s_1,\ell_{1,1}$ and $X=P$. The four properties are trivially fulfilled.

    For the induction step, we consider a caterpillar $G$ with a spine on $n\ge 2$ vertices. 
    Note that we can assume $s_n\geq 1$; otherwise we may shorten the spine.
    Let $G'$ denote the subcaterpillar obtained from $G$ by deleting the leaf~$\ell_{n,m_n}$. Observe that $G'$ is contained as a tube in our initial tubing, see \cref{fig:catB}. Let $Q'$ be the path we get by induction on $G'$ with the parts $Q_0'$ and $X'$. Each tubing in $Q'$ is maximally nested and hence associated to a permutation. Adding $\ell_{n,m_n}$ to the end of these permutations lifts the path $Q'$ to a path $Q$ in $\mathcal{A}(G)$. Along $Q$, all tubes not containing $s_n$ are dealt with, except the singletons that are leaves which we deliberately ignore for now. 
    Along $X'$, we see all tubes containing $s_n$ but not $\ell_{n,m_n}$. Again, writing $\ell_{n,m_n}$ at the end of each of the associated permutations lifts $X'$ to a path $\Hat{X}$ in $\mathcal{A}(G)$. Furthermore, let $\Bar{X}$ be the path $X'$ where we put $\ell_{n,m_n}$ into the second position of each permutation. Lastly, let $A$ be the path obtained by absorbing $\ell_{n,m_n}$ but skipping the last flip, i.e., not touching the singleton $s_n$, starting from the endpoint of $Q$.
    Then we define
    \[
        P = QA\Bar{X}\inv.
    \]
    To show that this is the required path, we first note that these paths fit together. By definition, the end of $Q$ and the start of $A$ match. The end of $A$ is the end of $X'$ with $\ell_{n,m_n}$ added in the second position of the permutation. This is exactly the starting point of $\Bar{X}\inv$.
    It is clear that the path
    \[
        X=(\Hat{X})A\Bar{X}\inv
    \]
    fulfills the conditions required. To show that $P_t$ is really connected for all tubes $t$ of $G$ that do not consist of only one leave of $G$, note that the tubes not containing $\ell_{n,m_n}$ are seen along $Q$ and those containing~$\ell_{n,m_n}$ are seen along $\Bar{X}$. To see that all tubes are seen we note that the tubes containing~$\ell_{n,m_n}$ (but not just containing $\ell_{n,m_n}$) are in bijection with those containing $s_n$ and not $\ell_{n,m_n}$. 

    To obtain a \fh path, start in the tubing that contains all leaves of $G$ as singletons as well as the tubes $T_i$ for $i=1,\ldots,n-1$, where $T_i$ contains all spine vertices $s_j$ with $j\le i$ and all leaves of the spine vertices it contains. Now flip to the starting point of $P$ in the obvious way and then apply $P$. This gives a \fh path. 
\end{proof}

\begin{figure}[hbt]
    \centering
    \includegraphics[scale = 0.32]{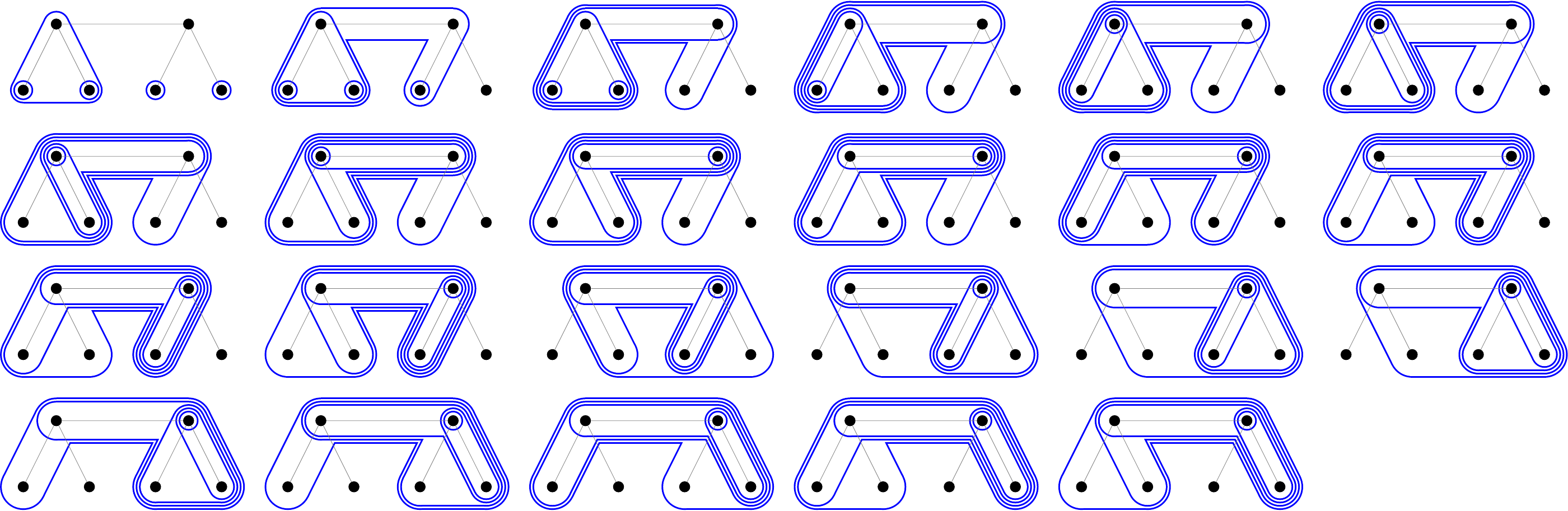}
    \caption{The \fh path constructed in the proof of \Cref{lem:caterpillar}. Note that this particular path can be closed to a cycle. It is not known whether this is always possible with this construction.}
    \label{fig:caterpillar}
\end{figure}

\cref{fig:caterpillar} illustrates a \fh path resulting from the constructive proof of \Cref{lem:caterpillar}.

\subsection{Complete bipartite graph associahedra} \label{sec:completeBipatite}

In this section, we prove the first part of \cref{thm:bipcat}, namely the following statement.

\begin{lemma}
    Graph associahedra of complete bipartite graphs have a \fh path.
\end{lemma}

We introduce some concepts.

\paragraph{Pattern.}
Let $G=K_{n,m}$ be a complete bipartite graph where the parts of the bipartition are $A=\{a_1,\dots,a_n\}$ and $B=\{b_1,\dots, b_m\}$. Consider a nested tubing in $\mathcal A(G)$ with permutation $\pi$. The \defn{pattern} of $\pi$ is defined to be the word $w\in \{A,B\}^{n+m}$ such that
\[
    w(i) = \begin{cases}
                A,& \text{ if }\pi(i)\in A\\
                B,& \text{ if }\pi(i)\in B
            \end{cases}.
\]

\paragraph{Shaving.}
Consider the nested tubing associated with the permutation \[a_1,b_1, a_2, \ldots, a_n,b_2, \ldots, b_m,\] with the pattern $ABA\ldots AB\ldots B$. Let $t$ denote the 2-tube $t=\{a_1,b_1\}$. The tubes containing $t$ are made of $t$ together with any subset of the other vertices. Hence, we have the structure of the permutahedron. Let the active set $M=A\setminus t$ consist of the elements of $A$ that are not in $t$. Apply the reverse path of the permutahedron $\mathcal A(K_{|M|})$ from \cref{thm:perm} to the active set, by modifying the permutation on these elements. Now take the first element of $B$ in the permutation that is not in~$t$ and flip it to the second place of the permutation (absorbing it). Add the replaced $B$ element to the active set. Apply the path of the permutahedron $\mathcal A(K_{|M|})$ to the now larger active set, again shifting in either direction. Now take the next $B$ element and flip it to the front. Add the replaced $B$ element to the active set and repeat. Do this until no new $B$ element can be flipped to the front. See \cref{fig:shaving} for an example. Finally, it is easy to check that if we always shift to the left the path ends in $a_1,b_n,a_3,\dots, a_n,a_2,b_1,\dots,b_{n-1}$.

\begin{figure}[htb]
    \centering
    \begin{align*}
        & & &a_1,b_1,\textcolor{blue}{a_2,a_3,a_4},b_2,b_3,b_4 & &\to &
        &a_1,b_1,\textcolor{blue}{a_3,a_4,a_2},\textcolor{red}{b_2},b_3,b_4\\
        &\to & &a_1,\textcolor{red}{b_2},\textcolor{blue}{b_1,a_3,a_4,a_2},b_3,b_4 & &\to &
        &a_1,b_2,\textcolor{blue}{a_3,a_4,a_2,b_1},\textcolor{red}{b_3},b_4\\
        &\to & &a_1,\textcolor{red}{b_3},\textcolor{blue}{b_2,a_3,a_4,a_2,b_1},b_4 & &\to &
        &a_1,b_3,\textcolor{blue}{a_3,a_4,a_2,b_1,b_2},\textcolor{red}{b_4}\\
        &\to & &a_1,\textcolor{red}{b_4},\textcolor{blue}{b_3,a_3,a_4,a_2,b_1,b_2} & &\to &
        &a_1,b_4,\textcolor{blue}{a_3,a_4,a_2,b_1,b_2,b_3}
    \end{align*}
    \caption{The sequence of permutations coming from the shaving procedure on $K_{4,4}$. The active set is drawn in blue and the element from $B$ that is absorbed is marked in red in each step. Note how the active set is always shifted to the left. The direction of shift can be chosen freely in each step.}
    \label{fig:shaving}
\end{figure}

\begin{lemma}
    The shaving procedure gives a valid path in $\mathcal A(K_{n,m})$ starting in the nested tubing associated with $a_1,b_1, a_2, \ldots, a_n, b_2, \ldots, b_m$ that visits every tube containing $a_1$ exactly once. 
\end{lemma}
\begin{proof}
    Let $T_1$ be the set of tubes containing $a_1$. Except for the singleton, each one contains at least one element of $B$ that ensures the connectedness. For a tube $t\in T_1\setminus \{a_1\}$ let $\max_B(t)$ be the maximal element in $B$ that $t$ contains. Let $T_1^{m}$ be the tubes in $T_1$ such that $\max_B(t)=m$. Starting with the permutation $
    a_1, b_1, a_2, \ldots, a_n, b_2, \ldots, b_m
    $ (pattern $ABA\ldots AB\ldots B$), in the first step we see exactly the tubes in $T_1^{b_1}$. Then we flip $b_2$ to the front. In the next shift we see exactly the tubes in $T_1^{b_2}$ and so forth. The fact that these sets partition $T_1$ proves the claim. 
\end{proof}

Note that the statement holds for any sequence of \defn{left/right shifts} we choose, where a right shift means that we apply the path obtained in \cref{thm:perm}, and a left shift means that we apply the reverse of that path. This has the effect of shifting the elements of the active set either to the right or to the left, respectively.

We use the shaving procedure to construct \fh paths in the graph associahedra of complete bipartite graphs.

\begin{proof}[Proof of \Cref{thm:bipcat} -- complete bipartite graphs]
    We start with the nested tubing corresponding to $a_1$, $b_1$, $a_2,\ldots, a_n, b_2,\ldots,b_n$ and apply the shaving procedure always shifting to the left. We end up with a tubing with the pattern $ABA\ldots AB\ldots B$. Now we have seen all tubes containing $a_1$. Flip the 2-tube, then we have two singletons on the $A$ side: $(AA)BA\ldots AB\ldots B$. Now flip $a_1$ all the way to the right (expelling it), we get the pattern $ABA\ldots AB\ldots B\| A$ where the $\|$ stands for the fact that we ignore $a_1$ from here on out. We see the same pattern on the left side of the $\|$, hence we proceed by induction. Eventually all elements of $A$ are expelled to the right side of the $\|$, leaving us with $AB\ldots B \| A\ldots A$. On the left side we have the star left. Flip the singleton, this leaves us with the pattern $BAB\ldots B$ for which the same procedure can be applied (with the roles of $A$ and $B$ reversed). In each step we see exactly those tubes containing a particular element of $A$, starting with $a_1$.
\end{proof}

\begin{figure}[htb]
    \centering
    \includegraphics[scale = 0.3]{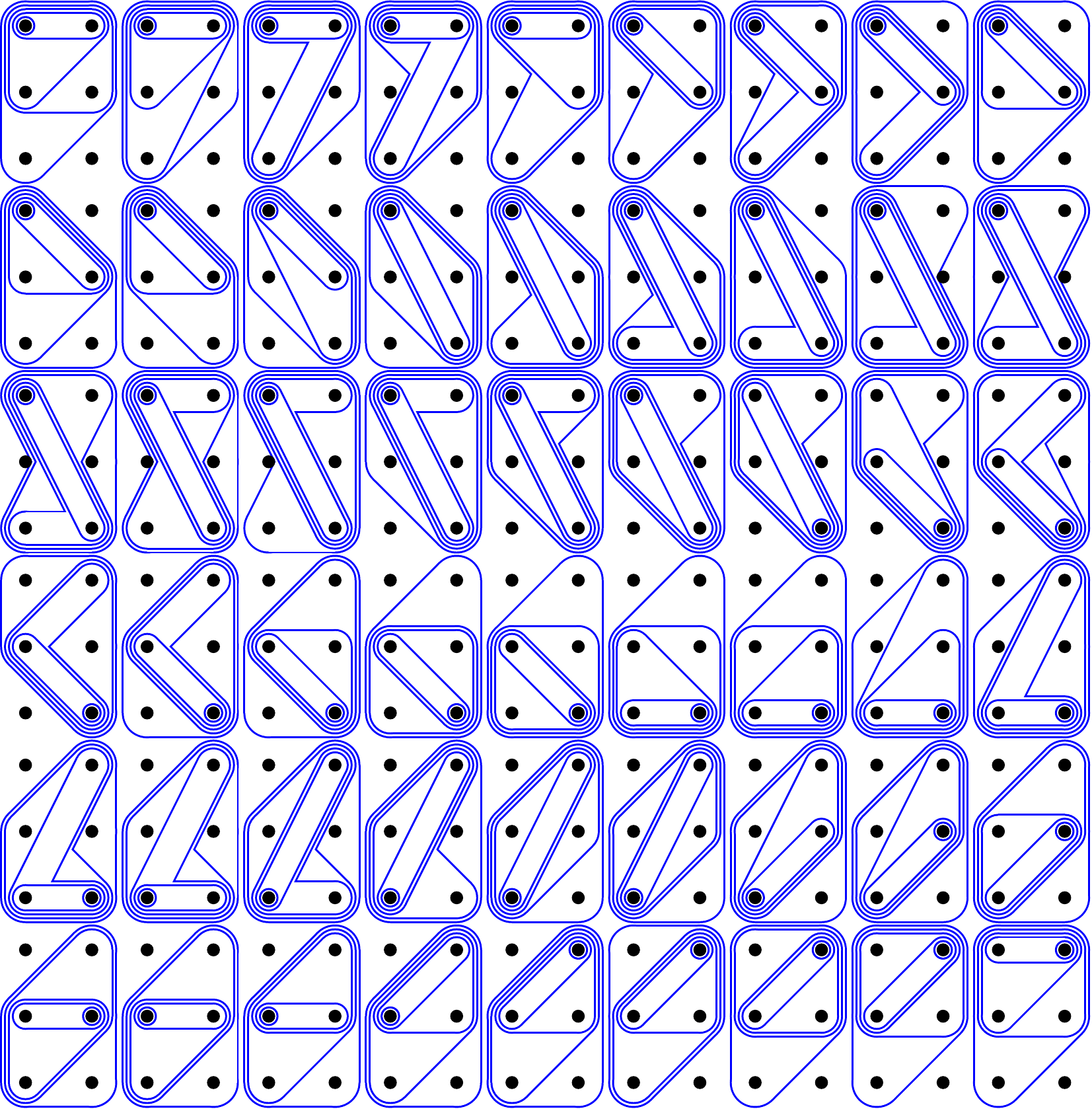}
    \caption{A \fh cycle in $\mathcal A(K_{3,3})$ (to be read from left to right in each line, the edges of $K_{3,3}$ are omitted for visibility reasons).}
    \label{fig:K33}
\end{figure}

\begin{figure}[htb]
    \centering
    \includegraphics[scale=0.8]{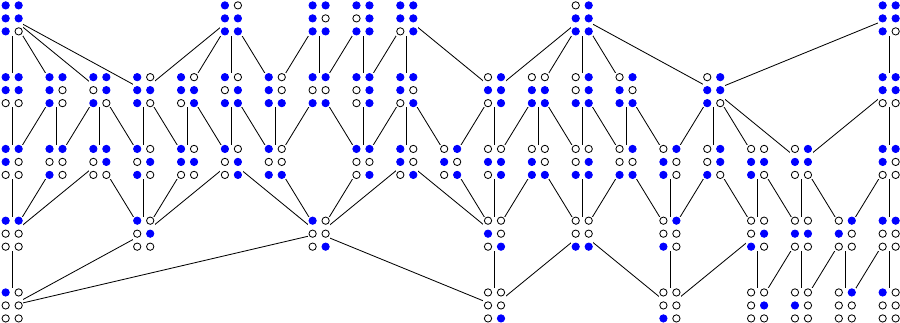}
    \caption{The rhombic strip for the cycle on $\mathcal A(K_{3,3})$ from \Cref{fig:K33}.}
    \label{fig:k33rhombic}
\end{figure}

\cref{fig:K33,fig:k33rhombic} illustrate a \fh cycle of $\mathcal A(K_{3,3})$ and its rhombic strip.


\section{On the existence of \fh paths and cycles} \label{sec:pathsCycles}

In the previous section, we presented some \fh paths. We now show that \fh simple 3-polytopes always admit a \fh path, but that the same is not true for non-simple polytopes.
\pathCycle*
\begin{proof}
We start by showing (i).
Consider a \fh cycle $C$ of $\mathcal P$ and let $k$ denote the number of facets.
A \fh cycle has length $k-3$ (and more generally $k-n$ for a $n$-polytope).
If there exists a facet that $C$ never leaves, then after the deletion of any two consecutive edges all faces are still guarded, and we obtain a  \fh path. Otherwise, by \cref{obs:simple}, 
$C$ has length~$k$ and we  aim to identify three edges, the deletion of which yields a \fh path.
In other words, we are seeking two consecutive vertices of $C$ such that their unique non-cycle edges lie on different sides of $C$: one inside and one outside $C$, as depicted in \cref{fig:cyclepath3DA}. When deleting the three corresponding edges of $C$, we obtain a path $P$ of length $k-3$ as desired. As the end vertices of~$P$ do not share a face (or equivalently, all facets are still guarded), $P$ is \fh. 
\begin{figure}[htb]
    \centering
    \begin{subfigure}[t]{0.3\textwidth}
    \centering
        \includegraphics[page=8]{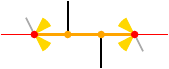}
        \subcaption{}
        \label{fig:cyclepath3DA}
    \end{subfigure}
    \hfil
    \begin{subfigure}[t]{0.3\textwidth}
    \centering
        \includegraphics[page=9]{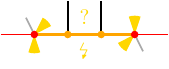}
        \subcaption{}
        \label{fig:cyclepath3DB}
    \end{subfigure}
    \hfil
    \begin{subfigure}[t]{0.3\textwidth}
    \centering
        \includegraphics[page=10]{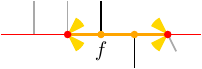}
        \subcaption{}
        \label{fig:cyclepath3DC}
    \end{subfigure}
    \caption{Illustration for the proof of \cref{{thm:pathCycle}}(ii). ({\sc a}) If the non-cycle edges of two consecutive vertices lie on different sides of $C$, then deleting their incident cycle edges yields a \fh path. ({\sc b}) If the non-cycle edges lie on the same side of $C$, then some face is not guarded. ({\sc c}) The last two vertices of $f$ have the desired property.}
    \label{fig:cyclepath3D}
\end{figure}

\cref{fig:cyclepath3DB} illustrates the case in which the two non-cycle edges lie on the same side of $C$. Then, in the corresponding path, one facet is not guarded and the end vertices lie on the same facet.

In this case, it remains to identify three suitable edges of $C$. For an illustration, consider \cref{fig:cyclepath3DC}. We consider some facet $f$ that is visited for at least two edges by $C$; it is easy to see that such a facet exists by 3-regularity. (In fact, by 3-regularity, for every edge of $C$, one of its two incident facets is visited at least by two edges of $C$.) Traversing $C$ in some direction, we consider the last two vertices that see~$f$. Clearly, the two incident non-cycle edges lie on different sides. This finishes the proof of~(i).

\medskip

It remains to show (ii).
Consider the 3-connected planar graph in \cref{fig:cycleVSpathB}. By Steinitz theorem, it is the skeleton of a 3-polytope $\mathcal P$. \cref{fig:cycleVSpathB} highlights a \fh cycle $C$. We remark that  this cycle is not unique; there also exist \fh  cycles of length 4.

\begin{figure}[htb]
    \centering
    \begin{subfigure}[t]{0.3\textwidth}
    \centering
        \includegraphics[page=2]{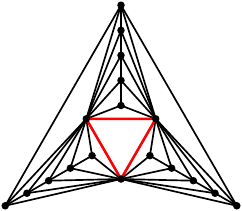}
        \subcaption{}
        \label{fig:cycleVSpathB}
    \end{subfigure}
    \hfil
    \begin{subfigure}[t]{0.15\textwidth}
    \centering
        \includegraphics[page=6]{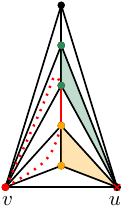}
        \subcaption{}
        \label{fig:cycleVSpathC}
    \end{subfigure}
    \hfil
    \begin{subfigure}[t]{0.15\textwidth}
    \centering
        \includegraphics[page=7]{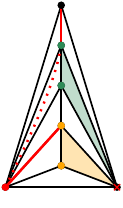}
        \subcaption{}
        \label{fig:cycleVSpathD}
    \end{subfigure}
    \hfil
    \begin{subfigure}[t]{0.15\textwidth}
    \centering
        \includegraphics[page=5]{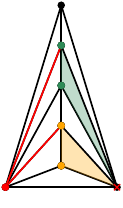}
        \subcaption{}
        \label{fig:cycleVSpathE}
    \end{subfigure}
    \hfil
    \caption{Illustration for the proof of \cref{{thm:pathCycle}}(ii).} 
    \label{fig:cycleVSpath}
\end{figure}

In the following we show that $\mathcal P$ does not have a \fh path.
Suppose that there exists a \fh path $P$ and let $\ell$ denote the number of vertices that $P$ and~$C$ have in common.
Clearly, $\ell\geq 1$; otherwise the central face is not visited.
If $\ell=3$, then exactly two  edges of $C$ belong to $P$; otherwise the central triangle is visited at least twice. However, then the triangle on the outer side of the unused edge 
is visited twice by $P$, a contradiction.

Thus $\ell\in\{1,2\}$, hence one vertex of $C$ is visited and one is not. Consider the subgraph illustrated in \cref{fig:cycleVSpathC} and suppose $P$ visits $v$ but not $u$. We aim to show that $v$ has two neighbors in $P$ that lie in this subgraph. In order to visit the green and orange faces,~$P$~visits an orange and a green vertex.
If~$P$ visits both after or before visiting $v$, then it uses the vertical (red) edge as illustrated in \cref{fig:cycleVSpathC}. However, as the incident faces must be visited consecutively, this yields a cycle (and degree two for~$v$ in this subgraph).  Therefore, assume some orange vertex is visited by an edge to $v$ and $P$ ends there, see \cref{fig:cycleVSpathD}.
Now, if $P$ visits no green vertex directly before or after $v$, then $P$ contains the topmost vertical (red) edge. However, as the incident faces are visited once, $P$ contains an edge to a green vertex. In particular, $v$ has two neighbors within the subgraph. (In fact, \cref{fig:cycleVSpathE} illustrates the unique solution for the subgraph, note that $P$ ends within the subgraph.)

If $\ell=2$, then $P$ contains the edge between these two vertices. However, by the above observation, the two vertices have degree two in their subgraph. A contradiction.
If $\ell=1$, then $v$ has degree two in each of the two subgraphs -- a contradiction again.
\end{proof}

\section{Computational complexity}
\label{sec:hardness}

\begin{figure}[ht]
	\centering
	\subfloat[\label{fig:vertex-breaka}]{\includegraphics[page=1]{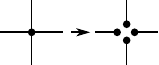}}\qquad
	\subfloat[\label{fig:vertex-breakb}]
  {\includegraphics[scale=.4]{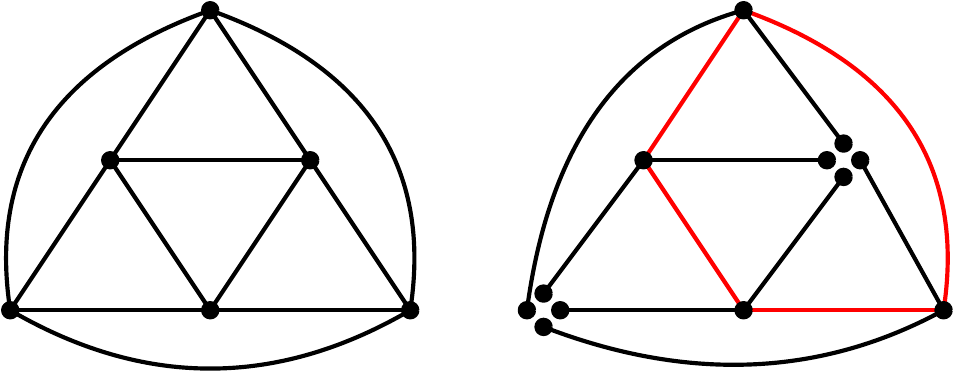}}
	\caption{({\sc a}) Vertex-breaking operation. ({\sc b}) Negative instance of TRVB.}
	\label{fig:vertex-break}
\end{figure}

In this section, we prove the \NP-completeness of the problem of deciding the existence of a \fh cycle in a given polytope.
We exhibit a reduction from the \NP-complete \textsc{Tree-Residue Vertex-Breaking} (TRVB) problem introduced by Demaine and Rudoy~\cite{demaine2018tree}. 
In TRVB, we are given a connected graph $G$ and asked whether we can transform $G$ into a tree by performing \defn{vertex-breaking} operations. 
A \defn{vertex-breaking} of a vertex $v$ replaces $v$ by $\deg(v)$ vertices of degree 1; for an illustration see \cref{fig:vertex-breaka}. 
The problem remains \NP-hard even if $G$ is a planar 4-regular graph.
\cref{fig:vertex-breakb} shows that the octahedral graph is a  negative instance of TRVB.
Note that if two adjacent vertices are broken, a component containing a single edge is created. Without loss of generality a potential solution breaks the leftmost vertex. Then, none of the vertices in the red cycle can be broken and we cannot eliminate all cycles of $G$ without disconnecting it.

\hardness*

\begin{proof}
	Containment in \NP\ is easy.
	\NP-hardness is established by a reduction from TRVB. 
	We first give an informal intuition. 
	Given a connected planar 4-regular graph $G'$ as instance of TRVB, we build a graph $G$ that will have a \fh cycle $C$ if and only if $G'$ is a positive instance of TRVB.
	We can think of $G$ as a ``thickening'' of $G'$, obtained by blowing up vertices and edges by their respective gadgets. 
	If $G'$ can be transformed into a tree $T$, the construction forces $C$ to act as an Euler tour around $T$ (using vertices of $G$).
	
	\begin{figure}[htb]
		\centering
		\subfloat[\label{fig:edge-gadget}]{\includegraphics[page=3]{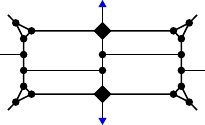}}\qquad
		\subfloat[\label{fig:cross-gadget}]{\includegraphics[page=4]{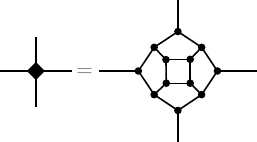} }\qquad
		\subfloat[\label{fig:vertex-gadget}]{\includegraphics[page=5]{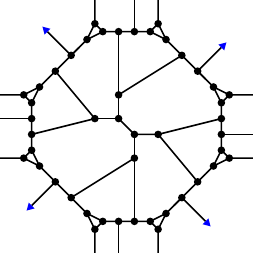}}
		\caption{({\sc a}) Edge gadget,  ({\sc b}) cross gadget , and ({\sc c}) vertex gadget.}
		\label{fig:gadgets}
	\end{figure}
 	
	The \defn{edge gadget} is shown in \cref{fig:edge-gadget}, and the \defn{vertex gadget} is shown in \cref{fig:vertex-gadget}.
	The leftmost (rightmost) 8 vertices of the edge gadget are shared with the vertex gadget corresponding to the left (right) endpoint of the edge.
	The diamond vertices in the edge gadget represent the subgraph shown in \cref{fig:vertex-gadget} that we call \defn{cross gadget}.
	The blue triangles denote edges whose other endpoints are not part of any gadget.
	Note that they point to an original face of $G'$.
	A face $F$ of $G'$ incident to $e_F$ edges has $2e_F$ associated blue triangles (one from every corresponding edge gadget and one for every corresponding vertex gadget). 
	Note that since $G'$ is simple, $e_F>2$.
	We connect the blue triangles of each face with a binary tree, similar to the 6 vertices in the interior of the vertex gadgets. 
	That means that $F$ correspond to $2e_F$ faces in $G$.
	The graph $G$ is planar and cubic.
    Note that two nonadjacent vertices are a 2-cut in a planar graph if and only if they both appear in two face cycles.
    By construction, all face cycles internal to gadgets intersect at exactly one edge. The only faces not internal to gadgets are the ones created by the addition of the binary trees connecting the blue triangles. By construction, such face cycles share exactly one edge with adjacent face cycles. We conclude that $G$ is 3-connected.
	This concludes the description of $G$.

It remains to show that $G$ is \fh if and only if $G'$ is a yes-instance of TRVB.
First, assume $G'$ is a yes-instance, so that it can be transformed into a tree $T$ by breaking a set of vertices.
	We now describe a corresponding \fh cycle $C$ in $G$.
	Informally, $C$ traces an Euler tour of $T$. 
	If a vertex $v$ of $G'$ was broken, we choose $C$ to traverse the corresponding vertex gadget as shown in Figure~\ref{fig:sol-break}.
	Otherwise, $C$ traverses the vertex gadget as shown in Figure~\ref{fig:sol-unbreak}.
	For every edge $e$ that had one of its endpoints broken, $C$ traverses the corresponding edge gadget as in Figure~\ref{fig:sol-edge1}. Else, none of the endpoints of $e$ are broken in $T$, and we make $C$ traverse the edge gadget as in  Figure~\ref{fig:sol-edge2}. All cross gadgets are traversed as in Figure~\ref{fig:sol-cross}.
	We can verify that the faces contained by the gadgets are traversed exactly once as desired. 
	The $2e_F$ faces corresponding to a face $F$ of $G'$ are all adjacent to a unique edge in an edge gadget that is incident to a cross gadget.
	Since $C$ contains those edges, these faces are visited. Each of these faces are incident to a single edge gadget and a single vertex gadget.
	By construction, $C$ also visits these faces only once in both cases: when the vertex is broken or not.	
	
	\begin{figure}[ht]
		\centering
		\subfloat[\label{fig:sol-break}]{
        \includegraphics[page=2]{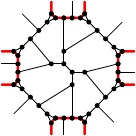}}\qquad
		\subfloat[\label{fig:sol-unbreak}]{
  \includegraphics[page=3]{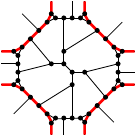}}\qquad
		\subfloat[\label{fig:sol-edge1}]{\includegraphics[page=8,scale=.7]{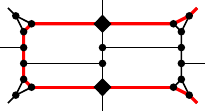}}\qquad
		\subfloat[\label{fig:sol-edge2}]{
  \includegraphics[page=9,scale=.7]{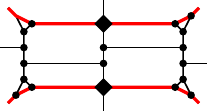}}\qquad
		\subfloat[\label{fig:sol-cross}]{
  \includegraphics[page=2]{figures/face-Ham-Hard2}}
		\caption{Constructing a \fh cycle given a solution to the instance of TRVB.}
		\label{fig:local-solution}
	\end{figure}

    For the reverse direction,  assume that $G$ admits a \fh cycle $C$.
	We now show that~$C$ must traverse~$G$ so that every vertex gadget is traversed exactly as in \cref{fig:sol-break} or \cref{fig:sol-unbreak}. That allows us to determine which vertices of $G'$ to break so that we have a solution for the TRVB problem.
	To this end, we prove some structural properties of $C$.

	\begin{figure}[ht]
		\centering
		\subfloat[\label{fig:cross-contradiction}]{
  \includegraphics[page=3]{figures/face-Ham-Hard2}}\qquad
		\subfloat[\label{fig:edge-contradiction}]{\includegraphics[page=12,scale=.7]{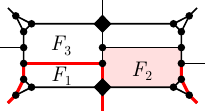}}\qquad
		\subfloat[\label{fig:vertex-contradiction}]{
   \includegraphics[page=4]{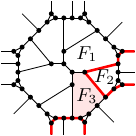}}\qquad
		\subfloat[\label{fig:edge-contradiction2}]{\includegraphics[page=13,scale=.7]{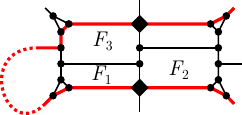}}\qquad
		\caption{Illustrations for the proofs of Claims~\ref{cl:cross}, \ref{cl:edge1} and \ref{cl:vertex}.
        }
		\label{fig:structure-C}
	\end{figure}
	 
	\begin{claim}\label{cl:cross}
		$C$  traverses exactly two of the four edges incident a cross gadget, and they must not be consecutive in the circular order around the gadget.
	\end{claim}
	\begin{proof}
		Because the cross gadget contains faces induced by the vertices of the gadget, it is clear that~$C$ must traverse at least two of the four incident edges or else these face would not be visited.
		For a contradiction, assume that $C$ traverses two consecutive edges.
		Refer to \cref{fig:cross-contradiction}.
		Without loss of generality, assume that $C$ traverses the left and top edges, which are both incident to face $F_1$. 
		Since the edges of $C$ on the boundary of $F_1$ must induce a path, only two faces enclosed by the gadget are visited. Then $C$ must visit the same gadget again, as shown in the figure. However $C$ cannot visit the small quadrangular face inside the gadget, a contradiction.
	\end{proof}

	\begin{claim}\label{cl:edge1}
		In every edge vertex, $C$  traverses the four horizontal (using the orientation of Figure~\ref{fig:edge-gadget}) edges incident to cross gadgets.
	\end{claim}
	\begin{proof}
		For contradiction, assume that $C$ traverses a vertical edge incident to a cross gadget. Refer to Figure~\ref{fig:edge-contradiction}.
		By Claim~\ref{cl:cross}, $C$ must traverse both vertical edges, one of each is incident to faces $F_1$ and~$F_2$ enclosed by the gadget.
		$C$ must also traverse the two small triangular faces at the bottom of the gadget, but it cannot use the two horizontal edges incident to the cross gadget.
		Since the traversal around $F_1$ must induce a path, $C$ must go through the three upper edges of $F_1$ as shown in the figure. However, $C$ must then visit $F_2$ twice, a contradiction.
	\end{proof}

	We now focus on vertex gadgets. Note that the vertex gadget contains a cycle with 32 vertices that enclose 6 vertices, 13 edges and 8 faces. We call these \defn{inner} vertices, edges and faces.

	\begin{claim}\label{cl:vertex}
		$C$ contains no inner edge of a vertex gadget.
	\end{claim}
	\begin{proof}
		To arrive at a contradiction, assume that $C$ contains an inner edge of a vertex gadget, and let~$P$ be a path in the gadget containing such an edge. Refer to \cref{fig:vertex-contradiction}.
		Note that $P$ must visit at least~3 inner faces (in the figure they are labeled $F_1$, $F_2$ and $F_3$). 
		Every inner face is adjacent to a small triangular face which is also present in the adjacent edge gadget. 
		By Claim~\ref{cl:edge1}, $P$ visits exactly two of these (at least three) triangular faces, because as soon as $P$ reaches the boundary of the vertex gadget it must continue to the adjacent edge gadget.
		But $C$ visits all of the triangular faces, implying that $C$ visits an inner face more than once ($F_3$ in the figure), a contradiction. 
	\end{proof}

	\begin{claim}\label{cl:edge}
		All vertex and edge gadgets are traversed as in \cref{fig:local-solution}. 
	\end{claim}
	\begin{proof}
		By Claim~\ref{cl:edge1}, the four horizontal edges incident to cross gadgets in an edge gadget must be in $C$.
		Claim~\ref{cl:vertex} rules out traversals such as the one shown in \cref{fig:edge-contradiction2}. 
		The horizontal edges in the middle of the gadget must also not be in $C$ or else either $F_1$ or $F_3$ is visited twice.
        Also note that only one edge of the triangular faces can be traversed: the bottom edge of $F_1$ is traversed by Claim~\ref{cl:edge1} and, if two edges of the neighbor triangular face is traversed, $F_1$ is visited twice.
		We are left with the traversals in \cref{fig:sol-edge1,fig:sol-edge2}.
		That implies (with Claim~\ref{cl:vertex}) that $C$ must be exactly as in Figures~\ref{fig:sol-break} 
        or~\ref{fig:sol-unbreak} at a vertex gadget. 
	\end{proof}

	By Claim~\ref{cl:edge} we can consistently decide which vertices to break based on $C$. Since $C$ is connected, every cycle in $G'$ has at least one vertex broken, and the set of broken vertices are not a cut set of~$G'$.
	Thus the result of breaking the selected vertices is a tree.
	That concludes the proof of the theorem.
\end{proof}

\subsection*{Acknowledgments} 
This work was initiated during the Tenth Annual Workshop on Geometry and Graphs, held at the Bellairs Research Institute in Holetown, Barbados in February 2023. The authors wish to thank the organizers and the participants, and in particular David Eppstein, Kolja Knauer, and Torsten Ueckerdt, for insightful discussions on this topic. They also thank the referees of a previous version of this manuscript, and in particular Vincent Pilaud, who pointed out the connection between \fh cycles and bipartite belts, and a referee of the current manuscript for pointing out the connection to Venn diagrams detailed in \Cref{subsubsec:venn}. 

\bibliographystyle{plainurl}
\bibliography{bib}

\end{document}